\documentstyle[epic,eepic,epsf]{article}

\newtheorem{example}{Example}[section]

\newtheorem{theorem}[example]{Theorem}
\newtheorem{corollary}[example]{Corollary}
\newtheorem{definition}[example]{Definition}

\newtheorem{lemma}[example]{Lemma}

\font\tenrm=cmr10
\def\cd{\cdots}

\def\La{\Lambda}

\def\P{{\cal P}}

\def\wt{\mbox{\sl wt}\,}

\def\Z{{\bf Z}}
\def\Zn{\Z_{\ge0}}

\def\FF{{\ssym F}}
\def\N{{\ssym N}}

\def\Z{{\ssym Z}}
\def\Q{{\ssym Q}}
\def\C{{\ssym C}}
\def\Proof{\medskip\noindent {\it Proof: }}

\def\adots{\mathinner{\mkern2mu\raise1pt\hbox{.} 
\mkern3mu\raise4pt\hbox{.}\mkern1mu\raise7pt\hbox{.}}}
\def\<{\langle}
\def\>{\rangle}

\def\SG{{\goth S}}

\def\mod{{\rm\,mod\,}}

\def\C{{\bf C}}

\def\ch{{\rm ch\,}}
\def\Sl{{\goth sl}}
\def\glchap{\widehat{\hbox{\goth gl}}}
\def\slchap{\widehat{\goth sl}}
\def\h{{\goth h}}

\def\F{{\cal F}}

\def\F{{\cal F}}
\def\dim{{\rm dim\,}}
\def\wt{{\rm wt\,}}

\def\F{{\cal F}}
\def\tr{{\rm tr\,}}

\def\boxit#1#2{\setbox1=\hbox{\kern#1{#2}\kern#1}%

\dimen1=\ht1 \advance\dimen1 by #1 \dimen2=\dp1 \advance\dimen2 by #1
\setbox1=\hbox{\vrule height\dimen1 depth\dimen2\box1\vrule}%
\setbox1=\vbox{\hrule\box1\hrule}%
\advance\dimen1 by .4pt \ht1=\dimen1
\advance\dimen2 by .4pt \dp1=\dimen2 \box1\relax}

\font\tensym=msbm10
\font\sevensym=msbm7
\font\fivesym=msbm5
\newfam\ssymfam
\textfont\ssymfam=\tensym
\scriptfont\ssymfam=\sevensym
\scriptscriptfont\ssymfam=\fivesym
\def\ssym{\fam\ssymfam\tensym}

\font\tengoth=eufm10
\font\sevengoth=eufm7
\font\fivegoth=eufm5

\newfam\gothfam
\textfont\gothfam=\tengoth
\scriptfont\gothfam=\sevengoth
\scriptscriptfont\gothfam=\fivegoth
\def\goth{\fam\gothfam\tengoth}

\font\twlrm=cmr12

\def\mathalign#1{\hbox to 0pt{\hss$\vcenter{\openup\jot
\tabskip=0pt plus1fil \halign to \hsize
{\hfil$\displaystyle ##$\tabskip=0pt&&$\displaystyle ##$\hfil
\tabskip=0pt plus1fil \cr #1}}$\hss}}

\def\makestrut#1#2{{\dimen12=#2
\divide\dimen12 by 4\dimen11=\dimen12\multiply\dimen11 by 3
\global\setbox#1=\hbox{\vrule height\dimen11 depth\dimen12 width0pt}}}

\newdimen\tadhdimen \newdimen\tabhdimen \newdimen\vdimen
\newdimen\smtadhdimen \newdimen\smtabhdimen
\newcount\splurt
\newbox\tadstrut \newbox\tabstrut
\newbox\smtadstrut \newbox\smtabstrut

\def\setyoungsize#1#2{          
  \tadhdimen=#1\tabhdimen=#1\advance\tabhdimen by -0.4truept%
  \vdimen=#2%
  \makestrut\tadstrut\vdimen
  \advance\vdimen by -0.4pt%
  \makestrut\tabstrut\vdimen}

\def\setsmyoungsize#1#2{        
  \smtadhdimen=#1\smtabhdimen=#1\advance\smtabhdimen by -0.4truept%
  \vdimen=#2%
  \makestrut\smtadstrut\vdimen
  \advance\vdimen by -0.4pt%
  \makestrut\smtabstrut\vdimen}

\setyoungsize{15pt}{15pt}
\setsmyoungsize{9pt}{9pt}

\def\youngt#1{%
  \vcenter{\offinterlineskip
  \halign{&\copy\tadstrut\hbox to \tadhdimen{\hss$##$\hss}\cr #1}}}
\def\youngd#1{%
  \vcenter{\offinterlineskip
  \halign{&\vrule##&\copy\tabstrut\hbox to \tabhdimen{\hss$##$\hss}\cr #1}}}

\def\smyoungt#1{{\vcenter{\offinterlineskip
  \halign{&\copy\smtadstrut
  \hbox to \smtadhdimen{\hss$\scriptstyle ##$\hss}\cr #1}}}}
\def\smyoungd#1{{\vcenter{\offinterlineskip
  \halign{&\vrule##&\copy\smtabstrut
  \hbox to \smtabhdimen{\hss$\scriptstyle ##$\hss}\cr #1}}}}

\def\hdashfill{\leaders\hbox to 3pt{%
\hfil\vrule width1.5pt height0.4pt depth0pt}\hfill}

\newcommand{\bi}{\begin{itemize}}
\newcommand{\ei}{  \end{itemize}}
\newcommand{\be}{\begin{equation}}
\newcommand{\ee}{  \end{equation}}
\newcommand{\bea}{\begin{eqnarray}}
\newcommand{\eea}{  \end{eqnarray}} 
\newcommand{\leftgroup }{\left  \lgroup}
\newcommand{\rightgroup}{\right \rgroup}

\newcommand{\infinity}{\infty}
\newcommand{\plusmin}{\pm}

\begin{document}

\title{Combinatorics of solvable lattice models, 
       and 
       modular representations of Hecke algebras
       }

\author{
O. Foda\thanks{
  Dept. of Mathematics, Univ. of Melbourne,
  Parkville, Victoria 3052, Australia.}, 
B. Leclerc\thanks{
  D\'ept. de Math\'ematiques, 
  Univ.  de Caen, 
  BP 5186, 14032 Caen cedex, France.}, 
M. Okado\thanks{
  Dept. of Math Sciences, 
  Fac. of Eng. Science, 
  Osaka U., 
  Osaka 560, Japan.}, 
J.-Y. Thibon\thanks{
  Inst. G. Monge, 
  Univ. de Marne-la-Vall\'ee,
  93166 Noisy-le-Grand cedex, France.}, 
T.A. Welsh$^{\ast}$
  }

\date{{\sl Dedicated to the memory of Abdus Salam}}

\maketitle

\begin{abstract}

We review and motivate recently-observed relationships between 
exactly solvable lattice models and modular representations of 
Hecke algebras. Firstly, we describe how the set of $n$-regular
partitions label both of the following classes of objects:

\noindent 1. The spectrum of unrestricted solid-on-solid lattice 
             models based on level-1 representations of the affine 
	     algebras $\slchap_n$,
           
\noindent 2. The irreducible representations of type-A Hecke 
             algebras at roots of unity: $H_m(\sqrt[n]{1})$.  

Secondly, we show that a certain subset of the $n$-regular partitions
label both of the following classes of objects: 

\noindent 1. The spectrum of restricted solid-on-solid lattice
             models based on cosets of affine algebras 
	     $C[n, 1, 1]$ $=$ 
             $(\slchap_n)_{1}$ $\times$
             $(\slchap_n)_{1}$ $/$
             $(\slchap_n)_{2}$,

\noindent 2. Jantzen-Seitz (JS) representations of 
             $H_m(\sqrt[n]{1})$: irreducible representations 
	     that remain irreducible under restriction to 
	     $H_{m-1}(\sqrt[n]{1})$. 

Using the above relationships, we characterise the JS representations 
of $H_m(\sqrt[n]{1})$ and show that the generating series that
count them are branching functions of $\slchap_n$. 

\end{abstract}

\vfill
\eject

\setcounter{secnumdepth}{10}
\setcounter{section}{-1}


\section{Introduction}

Exactly-solvable models in statistical mechanics and quantum 
field theory\footnote{For an overview of exactly solvable 
models, see \cite{ISZBook}.} offer concrete physical realisations 
of abstract mathematical structures. One can study them either 
from a physical or from a mathematical point of view.

From a physical point of view, one can say that though these 
models are highly non-trivial (they contain an infinite number 
of interacting degrees of freedom) they are 
solvable\footnote{Solvable here means that one can---at least in 
principle---compute any physical quantity {\it exactly}. Exactly 
here means {\it completely}, in the sense of not being an 
approximation of anything in some perturbation expansion. 
However, {\it exactly} does not mean {\it rigorously}, since 
these computations typically make heavy use of certain 
assumptions. In other words, 
in principle, our computations are either exactly correct, or 
exactly wrong! The issue is usually settled by heavy numerical 
computations, or even comparisons with actual experimental results. 
However, none of this need to worry us. In this work, we are only 
interested in mathematical aspects of these models.}
because they are based on rich and consistent algebraic 
structures\footnote{For an introduction to the algebraic approach 
to exactly solvable lattice models, see \cite{JimboMiwaBook}.}. 
The most notable of these structures are the affine and Virasoro 
algebras\footnote{Standard references on affine and Virasoro 
algebras are \cite{KacBook,KacRainaBook}. A particularly readable 
introduction to the subject is \cite{KassBook}.}. Physicists are 
typically not very surprised by that. They are used to the fact 
that underlying all meaningful physical models there exist 
consistent mathematics, and they take it for granted that 
mathematics is useful to physics. 

From a mathematical point of view, one can say that these models 
offer concrete realisations of highly abstract algebraic structures.  
Contrary to the physicists, mathematicians tend to be surprised by 
the fact that their abstractions turn out to be relevant to 
physical problems. 

These relationships are important because they imply 
that one can establish a dictionary between the two pictures: 
the physical and the mathematical. For example, one can establish 
a dictionary between the representation theory of Lie algebras, their 
affine extensions, and their deformations, on the one hand, and the theory 
of exactly solvable lattice models on the other hand.

The above relationships have been repeatedly used as a guiding principle 
that enriched both subjects. For example, computations of short distance 
correlations in lattice models are based entirely on techniques developed 
in the context of the representation theory of highest weight modules 
\cite{DFJMN,JMMN}, while intuition gained from the physics of solvable 
lattice models has initiated the theory of quantum groups 
\cite{ChariPressleyBook}. 

This motivates us to look for other examples of correspondences 
between the two subjects: exactly solvable models, and representation 
theory.

In \cite{FLOTW} we have described such a correspondence  
between the $(\slchap_n)_1 \times (\slchap_n)_1 / (\slchap_n)_2$
coset models and representations of the symmetric group 
over a field of characteristic $n$ and of Hecke algebras
at an $n$th root of 1
\footnote{Hecke algebras play an important role in 
the theory of exactly solvable models. For a review, see \cite{MartinBook}. 
The connections introduced in \cite{FLOTW} are specifically related to 
the theory of modular representations, such as the JS irreducibility 
conditions, the concepts of cores, etc.}. 
In this work, we present an introduction to \cite{FLOTW}, and 
a review of its context. 

\bigskip
Let us now describe briefly the contents of this paper.
In Section 1,
we introduce a class of exactly solvable 2-dimensional
statistical mechanical models: the ABF models, outline how
Baxter's corner transfer matrix (CTM) method can be used to
reduce the computation of certain physical quantities, the
1-point functions, to a problem in $q$-counting: computing the
generating function of a class of combinatorial objects called 
{\it paths}\footnote{Combinatorialists are familiar with these
objects. They are called {\sl Dyck paths}. However, the paths
that we are interested in counting will satisfy extra conditions.
In particular, they are required to start at a certain height,
end at a certain height, do not go below the $x$-axis, and do not
go above a certain model-dependent height.}. 

Solving the $q$-counting problem, we observe that the result is  
a character of the Virasoro algebra. This observation is by now
classical, and is due to \cite{DJMO1}. 
By the Goddard-Kent-Olive construction \cite{GKO},
the same expressions can
also be obtained as branching functions of tensor products of
level $\ell$ and level 1 representations of $\slchap_2$. 
It is in this sense that these
models can be regarded as coset models of type 
$C[2, \ell, 1] := (\slchap_2)_\ell \times (\slchap_2)_1 / (\slchap_2)_{\ell+1}$.

Once the connection with the affine algebra $\slchap_2$ is
observed, it is natural to seek the corresponding models that are
related to $\slchap_n$. 
In Section 2,
we introduce the representation theory and the combinatorics
underlying these models, the
JMO models based on the cosets $C[n, \ell, 1]$. 
In this case, the CTM method can also be applied to reduce the computation 
of 1-point functions to the $q$-enumeration of a class of weighted paths
\cite{DJMO1}. 
The relevant highest weight modules are now those 
of the ${\cal W}_n$-algebras whose characters arise as branching functions
of the coset $C[n, \ell, 1]$. 

In the case $\ell = 1$, another combinatorial description has been 
obtained in \cite{FOW} in terms of coloured partitions.
These partitions are characterised by certain simple conditions. 
Let us refer to these conditions as the FOW-conditions, and denote the
set of these partitions by $FOW(n, j, k)$, 
where $n$ stands for $\slchap_n$ and $\{j, k\}$ determine the 
branching function which is computed, namely 
$b_{\Lambda_j,\Lambda_0}^{\Lambda_k+\Lambda_{j-k}}$. 

It turns out that the FOW-partitions have a definite meaning in
the representation theory of the symmetric groups over a field
of characteristic $n$. In
particular, if one disregards the $\{j, k\}$ labels in the
FOW-conditions, one obtains the conditions satisfied by the
Jantzen-Seitz (JS) partitions \cite{JS}. These are the partitions
that label irreducible representations of a symmetric group 
over $\FF_n$ that remain irreducible under restriction to
the symmetric subgroup of immediately lower rank.
Let us denote these partitions by $JS(n)$.
 
The first goal of this work is to explain why 
\begin{equation}\label{GOAL}
\coprod_{j, k} FOW(n, j, k) = JS(n).
\end{equation}

In Section 3,
we review some combinatorial aspects of the representations 
of symmetric groups over fields of characteristic $n$
and Hecke algebras $H_m(\sqrt[n]{1})$. 
Here of course, $n$ is assumed to be prime in the case
of symmetric groups.
Because our statistical mechanical
models are based on affine algebras $\slchap_n$ for all
$n>1$, the Hecke algebras form the natural context in which we
can work.
We describe in particular the Jantzen-Seitz problem for representations 
of symmetric groups, whose complete solution was
recently given by Kleshchev \cite{Kl1}, and formulate a similar problem for
Hecke algebras.

One of the main ingredients in our solution of the Hecke JS
problem is a relationship between the representation theory 
of Hecke algebras $H_m(\sqrt[n]{1})$
and the basic representation of the affine algebra 
$\slchap_n$. 
Namely, the sum over $m$ of the Grothendieck groups
associated with $H_m(\sqrt[n]{1})$ can be identified
with the irreducible representation $V(\Lambda_0)$ of $\slchap_n$, 
with Robinson's $i$-induction and $i$-restriction operators acting
as the Chevalley generators.
This relationship was described in \cite{LLT1,LLT}. 
Moreover, it was conjectured there that Kashiwara's upper
global crystal basis of $V(\Lambda_0)$ coincides
with the natural basis of the sum of Grothendieck groups
given by the classes of simple modules. This conjecture was proven 
by Ariki \cite{Ar2}, and independently by Grojnowski
using results of \cite{Gr}. 

In Section 4,
 we review the basic facts that we need from the
representation theory of quantum affine algebras. In particular, we
describe the Fock space representations of
$U_q(\slchap_n)$.
The theory of crystal bases is reviewed in
Section 5, 
and its relevance to the representation
theory of Hecke algebras is explained.

In Section 6, using the above-mentioned relationship,
we show that the Jantzen-Seitz problem for the 
Hecke algebras $H_m(\sqrt[n]{1})$
is equivalent to the computation via crystal graphs of the decomposition 
of tensor products of level 1 $\slchap_n$-modules.
These are precisely the tensor products arising 
in the coset models $C[n, 1, 1]$ studied in \cite{FOW}.
Therefore we obtain the desired explanation of (\ref{GOAL}). 
Moreover, building upon results of \cite{FOW}, we can express  
the generating functions of JS-partitions having a given
$n$-core in terms of branching functions of $\slchap_n$.
This is to be compared with a well-known result on
blocks of Hecke algebras. 
Indeed, the blocks of $H_m(\sqrt[n]{1})$ are
labelled by $n$-cores, and the dimension
of a block is the number of $n$-regular partitions of $m$
with the corresponding $n$-core. 
Using a formula first proved by Robinson (in the symmetric
group case) one can compute the generating series
of the dimensions of all blocks labelled by a given
$n$-core, and recognize the string function of the level 1
$\slchap_n$-modules.

Thus our result shows that some branching functions of $\slchap_n$ other
than the level 1 string function arise in a natural way in the modular
representations of $H_m(\sqrt[n]{1})$.
We conclude by discussing possible generalizations (Section 7).

\medskip
In Appendix A, we recall a description of the Specht modules of 
$H_m(v)$.

\section{From statistical mechanics to 
         Virasoro highest weight modules}

\subsection{Introduction}

A central activity in statistical mechanics is to describe critical 
phenomena in terms of exactly solvable lattice models. The most 
well-known of these models is the 2-dimensional Ising model, whose 
partition function was exactly calculated by Onsager in 1944 \cite{O}. 

In 1980, the hard hexagon model was solved by Baxter \cite{B}. In 
calculating order parameters, he found remarkable connections with 
pure mathematics: $q$-series identities of the Rogers-Ramanujan
type, modular functions, etc. In 1984, this model was extended by 
Andrews, Baxter and Forrester (ABF) to an infinite series of models
\cite{ABF}.

What ABF did can be summarised as follows. They constructed a series 
of restricted solid-on-solid (RSOS) models labelled by a positive 
integer $L$. Their Boltzmann weights are parametrised using elliptic 
functions, and satisfy the star-triangle relation or Yang-Baxter 
equation. Using the corner transfer matrix (CTM) method, they 
calculated certain physical quantities called local height 
probabilities (LHP's). From the latter, they computed the critical 
exponents: physical quantities that can be measured in actual
experiments. 

Subsequently, the Kyoto group \cite{DJKMO} realised that the generating 
functions of certain combinatorial objects: the so-called {\it paths}, 
or {\it 1-dimensional configurations}, that appear in the course of 
calculating the LHP's are precisely the branching functions for
cosets of the affine algebras: 
$(\slchap_2)_l\times(\slchap_2)_1/(\slchap_2)_{l+1}$,
where the subscripts stand for the level of the representation. 

This connection with affine algebras brought us an exciting game to 
play \cite{DJMO1}. The rules of this game\footnote{Also known as
{\it Baxter's ball game.}} are as follows\footnote{Although we only 
consider $(\slchap_n)$ case here, the situations for the other 
affine Lie algebras are similar. There is also a vertex model 
version of this game. In this case, the coset should be replaced 
by the level $l$ representation of the affine algebra $(\slchap_n)$,
and the branching function by a string function.}: 

\begin{itemize}
\item[(R1)] Find a solution of the Yang-Baxter equation related
            to the affine coset 
$(\slchap_n)_{l_1}\times(\slchap_n)_{l_2}/(\slchap_n)_{l_1+l_2}$.

\item[(R2)] Show that the 1-dimensional configuration sum is  
            a branching function for this coset.
\end{itemize}

In this work, we will not delve into the details of the above two
steps. We will rather consider the solutions of the Yang-Baxter 
equations as given, outline the way that the CTM produces
1-dimensional configurations in the simpler case of the ABF
models, and go straight to the combinatorial objects that it
produces. 

\subsection{The ABF restricted solid-on-solid models} 

A restricted solid-on-solid (RSOS) model, as defined in \cite{ABF}, 
is a system of interacting discrete variables, called \lq heights\rq, 
$l$, defined on the sites of a square lattice. The heights, 
$l$, take values in the set of positive integers 
$l \in \{1, 2, \cdots, L - 1\}$. $L$ characterises the model. 

RSOS models are interaction-round-face (IRF) models: each   
configuration of 4 heights defined on the corners of a face on the 
lattice is assigned a Boltzmann weight. In the notation of \cite{ABF}, 
the non-vanishing weights are
\bea
W(l, l + 1\vert l - 1, l) &=& W(l, l - 1\vert l + 1, l) = \alpha_l \nonumber \\
W(l + 1, l\vert l, l - 1) &=& W(l - 1, l\vert l, l + 1) = \beta_l  \nonumber \\
W(l + 1, l\vert l, l + 1) &=& \gamma_l                             \nonumber \\
W(l - 1, l\vert l, l - 1) &=& \delta_l                          
\eea
the labelling being given by
\begin{center}
\unitlength=0.00625in
\begin{picture}(254,141)(0,-10)
\path(20,120)(120,120)(120,20) (20,20)(20,120)
\put(0,120){\makebox(0,0)[lb]{\raisebox{0pt}[0pt][0pt]{\shortstack[l]{{\tenrm m}}}}}
\put(140,120){\makebox(0,0)[lb]{\raisebox{0pt}[0pt][0pt]{\shortstack[l]{{\tenrm n}}}}}
\put(0,0){\makebox(0,0)[lb]{\raisebox{0pt}[0pt][0pt]{\shortstack[l]{{\tenrm i}}}}}
\put(140,0){\makebox(0,0)[lb]{\raisebox{0pt}[0pt][0pt]{\shortstack[l]{{\tenrm j}}}}}
\put(180,60){\makebox(0,0)[lb]{\raisebox{0pt}[0pt][0pt]{\shortstack[l]{{\tenrm
$= W(m, n \mid i, j)$}}}}}
\end{picture}
\end{center}
The Boltzmann weights are parametrised as
\bea
\alpha_l &=& \rho h ( v + \eta ) \nonumber \\
\beta_l  &=& \rho h ( v - \eta ) 
\frac{
\left\lgroup h( w_{l - 1} )  h( w_{l + 1} ) \right\rgroup^{1/2} 
}{
h( w_{l} )
}
\nonumber \\
\gamma_l &=& \rho h ( 2 \eta )  
\frac{
 h( w_{l} + \eta - v )
}{
h( w_{l} )
} 
\nonumber \\  
\delta_l &=& \rho h ( 2 \eta )  
\frac{
 h( w_{l} - \eta + v )
}{
h( w_{l} )
}
\eea
\bea
h(u) = 2 p^{1/4} \sin {\pi u} \prod_{n = 1}^{\infinity}
\left\lgroup 1 - 2 p^n \cos \pi u + p^{2 n} \right\rgroup
\left\lgroup 1 - p^{2 n} \right\rgroup^{2}
\eea
$\rho$ is a normalisation factor, $w_l = 2 l \eta$.
They satisfy the star-triangle, or Yang-Baxter equation
\bea
\sum_{g}
W(b, c\, \vert a, g) W^{\prime}(a, g\, \vert f, e) W^{\prime \prime}(g, c\, \vert e, d)
\nonumber \\
=
\sum_{g}
W^{\prime \prime}(a, b\, \vert f, g) W^{\prime}(b, c\, \vert g, d) W(g, d\, \vert f, e)
\eea
They are characterised by the property that the heights on 
neighbouring sites differ by $\plusmin 1$. Because of this, 
one obtains a non-trivial model
only for $L - 1\geq 3$. For $L - 1 = 3$, one obtains the IRF formulation of 
the Ising model. In the 
limit $L \rightarrow \infinity$ one recovers the solid-on-solid version
of Baxter's 8-vertex model \cite{BaxterBook}. For a complete discussion of the RSOS 
models, we refer to \cite{ABF}.

\subsection{Corner transfer matrices}

One can generate configurations in a model defined on a square lattice 
using corner transfer matrices CTM's\footnote{There are two basic
references to the CTM methods: 1. Chapter 13 of \cite{BaxterBook}, and 
Appendix A of \cite{ABF}. The first gets more into the physical
background and assumptions that underlie the method. The second
is a shorter, more technical outline of the method.}. 
One divides the square lattice into 4 quadrants. As we want to 
compute certain quantities in the infinite lattice limit, we may 
assume that each quadrant of the finite lattice is a right triangle 
bounded by two segments of the axes, with the external spins fixed 
to some ground state boundary condition, as in the following picture
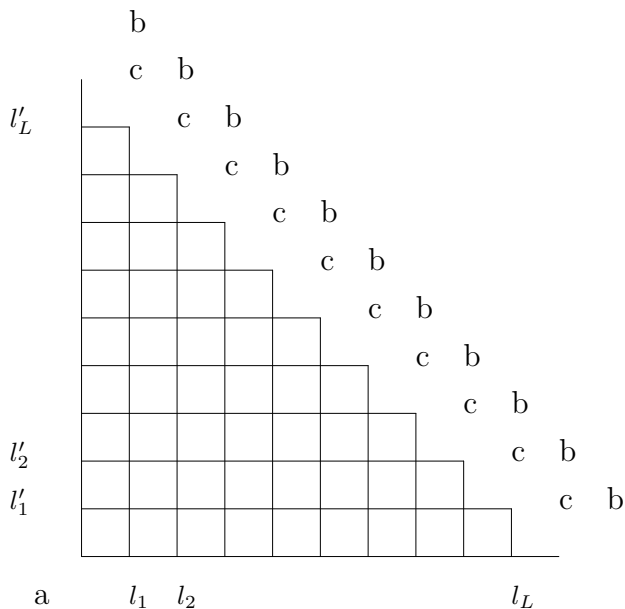
\begin{center}
\unitlength=0.0125in
\begin{picture}(256,264)(0,-10)
\path(30, 220)( 30,20)
\path(50, 200)( 50,20)
\path(70, 180)( 70,20)
\path(90, 160)( 90,20)
\path(110,140)(110,20)
\path(130,120)(130,20)
\path(150,100)(150,20)
\path(170, 80)(170,20)
\path(190, 60)(190,20)
\path(210, 40)(210,20)
\path(30,200)(50,200)
\path(30,180)(70,180)
\path(30,160)(90,160)
\path(30,140)(110,140)
\path(30,120)(130,120)
\path(30,100)(150,100)
\path(30,80)(170,80)
\path(30,60)(190,60)
\path(30,40)(210,40)
\path(30,20)(230,20)
\put( 50,240){\makebox(0,0)[lb]{\raisebox{0pt}[0pt][0pt]{\shortstack[l]{{\twlrm b}}}}}
\put( 70,220){\makebox(0,0)[lb]{\raisebox{0pt}[0pt][0pt]{\shortstack[l]{{\twlrm b}}}}}
\put( 90,200){\makebox(0,0)[lb]{\raisebox{0pt}[0pt][0pt]{\shortstack[l]{{\twlrm b}}}}}
\put(110,180){\makebox(0,0)[lb]{\raisebox{0pt}[0pt][0pt]{\shortstack[l]{{\twlrm b}}}}}
\put(130,160){\makebox(0,0)[lb]{\raisebox{0pt}[0pt][0pt]{\shortstack[l]{{\twlrm b}}}}}
\put(150,140){\makebox(0,0)[lb]{\raisebox{0pt}[0pt][0pt]{\shortstack[l]{{\twlrm b}}}}}
\put(170,120){\makebox(0,0)[lb]{\raisebox{0pt}[0pt][0pt]{\shortstack[l]{{\twlrm b}}}}}
\put(190,100){\makebox(0,0)[lb]{\raisebox{0pt}[0pt][0pt]{\shortstack[l]{{\twlrm b}}}}}
\put(210, 80){\makebox(0,0)[lb]{\raisebox{0pt}[0pt][0pt]{\shortstack[l]{{\twlrm b}}}}}
\put(230, 60){\makebox(0,0)[lb]{\raisebox{0pt}[0pt][0pt]{\shortstack[l]{{\twlrm b}}}}}
\put(250, 40){\makebox(0,0)[lb]{\raisebox{0pt}[0pt][0pt]{\shortstack[l]{{\twlrm b}}}}}
\put( 50,220){\makebox(0,0)[lb]{\raisebox{0pt}[0pt][0pt]{\shortstack[l]{{\twlrm c}}}}}
\put( 70,200){\makebox(0,0)[lb]{\raisebox{0pt}[0pt][0pt]{\shortstack[l]{{\twlrm c}}}}}
\put( 90,180){\makebox(0,0)[lb]{\raisebox{0pt}[0pt][0pt]{\shortstack[l]{{\twlrm c}}}}}
\put(110,160){\makebox(0,0)[lb]{\raisebox{0pt}[0pt][0pt]{\shortstack[l]{{\twlrm c}}}}}
\put(130,140){\makebox(0,0)[lb]{\raisebox{0pt}[0pt][0pt]{\shortstack[l]{{\twlrm c}}}}}
\put(150,120){\makebox(0,0)[lb]{\raisebox{0pt}[0pt][0pt]{\shortstack[l]{{\twlrm c}}}}}
\put(170,100){\makebox(0,0)[lb]{\raisebox{0pt}[0pt][0pt]{\shortstack[l]{{\twlrm c}}}}}
\put(190, 80){\makebox(0,0)[lb]{\raisebox{0pt}[0pt][0pt]{\shortstack[l]{{\twlrm c}}}}}
\put(210, 60){\makebox(0,0)[lb]{\raisebox{0pt}[0pt][0pt]{\shortstack[l]{{\twlrm c}}}}}
\put(230, 40){\makebox(0,0)[lb]{\raisebox{0pt}[0pt][0pt]{\shortstack[l]{{\twlrm c}}}}}
\put( 10,  0){\makebox(0,0)[lb]{\raisebox{0pt}[0pt][0pt]{\shortstack[l]{{\twlrm a}}}}}
\put( 50,  0){\makebox(0,0)[lb]{\raisebox{0pt}[0pt][0pt]{\shortstack[l]{{\twlrm $l_1$}}}}}
\put( 70,  0){\makebox(0,0)[lb]{\raisebox{0pt}[0pt][0pt]{\shortstack[l]{{\twlrm $l_2$}}}}}
\put(210,  0){\makebox(0,0)[lb]{\raisebox{0pt}[0pt][0pt]{\shortstack[l]{{\twlrm $l_L$}}}}}
\put( 0,  40){\makebox(0,0)[lb]{\raisebox{0pt}[0pt][0pt]{\shortstack[l]{{\twlrm $l'_1$}}}}}
\put( 0,  60){\makebox(0,0)[lb]{\raisebox{0pt}[0pt][0pt]{\shortstack[l]{{\twlrm $l'_2$}}}}}
\put( 0, 200){\makebox(0,0)[lb]{\raisebox{0pt}[0pt][0pt]{\shortstack[l]{{\twlrm $l'_L$}}}}}
\end{picture}
\end{center}
Starting from a half-axis that 
extends from the centre of the lattice to a corner, {\it e.g.} the 
south one, a corner transfer matrix generates height configurations 
on the quadrant bounded by the initial half-axis and the half-axis 
that follows, as one turns the initial half-axis in a clockwise 
direction, keeping one end fixed at the centre of the lattice.

We denote the heights on the initial half-axis by the {\it height 
vector} ${\bf l} = \{l_1, l_2, \cdots, l_m \}$, and those on the 
final half-axis by  
${\bf l^{\prime}} = \{l_{1}^{\prime}, l_2^{\prime}, \cdots, l_m^{\prime}\}$, 
with $l_1 = l_1^{\prime}$, where $l_1$ is the height at the centre of 
the lattice.

The elements of the corner transfer matrix are sums over products of 
Boltzmann weights for configurations labelled by fixed initial and final 
height vectors. More precisely, the matrix elements
of the corner transfer matrix  $A$ of the lower right quadrant
are defined by
$$
A_{{\bf l},{\bf l^{\prime}}}
=A_{{\bf l},{\bf l^{\prime}}}(b,c)
=
\left\{\matrix{
\sum\prod W(\sigma_i\sigma_j|\sigma_m\sigma_n)\quad {\rm if}\ \ l_1=l'_1 \cr
0 \ \ {\rm otherwise}
}
\right.
$$
where the product is over all faces of the quadrant, and the
sum runs over all admissible values of the heights on the
circled vertices, the other boundary spins being fixed
to $b,c$ as indicated on the picture.

One can think of the elements of a corner transfer matrix 
as the  partition functions of a quadrant of a square lattice with
fixed boundary conditions.    

We denote the corner transfer matrices of the 4 quadrants by
$A, B, C, D$ (in counterclockwise order).  

We take $m \rightarrow \infinity$, in the thermodynamic
limit.

\subsection{Local height probabilities} 

One can use the corner transfer matrices to compute 1-point 
functions, or equivalently order parameters. These functions 
are also known as the local height probabilities (LHP's)
$P\left\lgroup a \vert b, c \right\rgroup$.

Following \cite{BaxterBook,ABF}, we consider a RSOS model on 
a square lattice and fix 
as above
the height variables at the boundary 
of the lattice to the ground-state pair $\{b, c\}$, where 
$b - c = \plusmin 1$. The partition function 
$Z\left\lgroup b, c \right\rgroup$ of this configuration is 
\bea
Z \left\lgroup b, c       \right\rgroup = {\rm Trace} 
  \left\lgroup A B C D    \right\rgroup   \label{normalization}
\eea
The product of the corner transfer matrices in (\ref{normalization}) 
indicates identifying the heights on the common half-axes and 
summing over all possible heights. The trace indicates doing the 
same thing for the initial half-axis of $A$, and the final 
half-axis of $D$.

Computing $Z\left\lgroup b, c \right\rgroup$, using CTM's, allows 
us to fix the height at the centre of the lattice. We can do that 
by inserting the operator 
\bea
{\left\lgroup S_a \right\rgroup}_{ {\bf l,l^{\prime} } } = 
\delta\left\lgroup     {l_1}, a\right\rgroup
\delta\left\lgroup {\bf l, l^{\prime} } \right\rgroup ,
\quad
\delta\left\lgroup {\bf l, l^{\prime} } \right\rgroup \equiv 
\prod_{k = 1}^{m} \delta \left\lgroup  l_k, l^{\prime}_k  \right\rgroup       
\label{operator} 
\eea 
in (\ref{normalization}) to obtain 
\bea
Z\left\lgroup a \vert b, c\right\rgroup =
{\rm Trace} \left\lgroup S_a A B C D \right\rgroup  
\label{partitionfunction}
\eea
$Z\left\lgroup a \vert b, c\right\rgroup$ is the partition function 
of the configuration with fixed boundary conditions both at the centre 
and the boundary. `Local height' refers to the height at the centre of 
the lattice. To obtain a probability, we normalise 
$Z\left\lgroup b, c
  \right\rgroup$ 
  by 
  $Z\left\lgroup a \vert b, c
    \right\rgroup$:
\bea
P \left \lgroup a \vert b, c \right \rgroup =
\frac{
{\rm Trace } \leftgroup S_{a} A B C D \rightgroup
     }{
{\rm Trace } \leftgroup       A B C D \rightgroup
      }                         \label{LHP}
\eea  
where $P\left\lgroup a \vert b, c \right\rgroup$ is the probability that 
the height variable at the centre of the lattice is $a$, while the heights 
on the boundary are $\{b,c\}$. 
It satisfies
\bea
\sum_{a}  P\left\lgroup a \vert b, c \right\rgroup = 1  \label{star1}
\eea

\subsection{From 2-dimensional to 1-dimensional configurations}

We wish go through the derivation of LHP's, following \cite{ABF}.
Using the Yang-Baxter equations one can show 
that\footnote{A better phrasing may be {\it "one can argue
that"}.}, in the thermodynamic 
limit (where the lattice becomes infinitely large), the CTM's can be 
written in the form 

\begin{eqnarray}
A(v) \sim \alpha(v)Q_{1} M_{1} & e^{         v {\cal H}} & Q_{2}^{-1} \nonumber \\
B(v) \sim \beta(v)Q_{2}  M_{2} & e^{     -   v {\cal H}} & Q_{3}^{-1} \nonumber \\  
C(v) \sim \gamma(v)Q_{3} M_{3} & e^{         v {\cal H}} & Q_{4}^{-1} \nonumber \\
D(v) \sim \delta(v)Q_{4} M_{4} & e^{     -   v {\cal H}} & Q_{1}^{-1} 
\label{star3}
\end{eqnarray}

\noindent where ${\cal H}, Q_{1}, \cdots, Q_{4}, M_{1}, \cdots, M_{4}$, 
are matrices that are independent of the spectral parameter $v$. Further,
$\alpha,\beta,\gamma,\delta$ are scalar functions,
and 
${\cal H}, M_{1}, \cdots, M_{4}$ are diagonal. Substituting in (\ref{LHP}), 
we obtain

\bea
P\left\lgroup a \vert b, c \right \rgroup 
    = 
\frac{ 
{\rm Trace} {\leftgroup S_a M_1 M_2 M_3 M_4 \rightgroup}
             }{ 
{\rm Trace} {\leftgroup     M_1 M_2 M_3 M_4 \rightgroup}                    
              }
\eea

Following a series of arguments given in detail in \cite{BaxterBook}, 
one can show that, up to irrelevant scalar factors, 

\bea
A(\eta)    &=& C(\eta) = I, \nonumber \\ 
B( - \eta) &=& D( - \eta) = R_1 \nonumber \\
\eea
where 
$$
\left\lgroup R_j \right\rgroup_{ {\bf l, l^{\prime} } } =
\left\lgroup h(w_{ l_j }) \right\rgroup^{1/2}  \delta( {\bf l, l^{\prime} } )\, ,
$$
and $\eta$ is the crossing parameter. Substituting the above in (\ref{star3}), 
and using the fact that both the $Q_i$'s and ${\cal H}$ are diagonal, we 
obtain\footnote{This entire discussion is restricted to the anti-ferromagnetic, 
or so-called regime-III of the ABF models.}
\bea
A(\eta) B(- \eta) C(\eta) D( - \eta) = R_1^2
\eea
\bea
M_1 M_2 M_3 M_4 = R_1^2 e^{- 4 \eta \cal{H} }. \label{MMMM}
\eea

Evaluating the matrices in (\ref{MMMM}), one obtains
\bea
& 
\left\lgroup e^{- 4 \eta {\cal H}}
\right\rgroup_{{\bf l l^{\prime} }}  = &
                        q^{ 
\left\lgroup - (2 l_1 - L)^2 / 16 L +  \Phi_m ({\bf l}) 
\right\rgroup} 
                        \delta ({\bf l, l^{\prime}})          
\eea
where
\bea
\Phi_{m}(l_1, \cdots, l_{m+2})
                & = & \sum_{j = 1}^{m}
                      j {\cal H}(l_j, l_{j +1}, l_{j +2})
                                                    \nonumber \\
{\cal H}(l_j, l_{j + 1}, l_{j + 2}) & = & \frac{1}{4} \vert l_j - l_{j + 2}\vert
                                                    \nonumber \\
\eea
and
\bea
{\left \lgroup R_{1}^{2 } \right \rgroup}_{l l^{\prime}} =
\tau q^{{( l_1 - L)}^{2}/ 8 L} E (q^{l_1
/ 2}, q^{L/2} ) 
\delta ({\bf l}, {\bf l^{\prime} })                                  
\eea
\noindent where $\tau = - \log q / (2 \pi)$, and  
$            
E(z, q) = \Pi_{n=1}^{\infinity}(1 - q^{n - 1}z)(1 - q^n z^{- 1})(1 - q^n).
$
Substituting the above results in (\ref{normalization}), and (\ref{LHP}), 
one obtains

\bea
& &
P \leftgroup a \vert b, c \rightgroup  = 
\frac{
E(q^{a/2}, q^{L/2}) X_m(a, b, c| q)
     }{
\sum_{1 \leq d < r} E(q^{d/2}, q^{L/2}) X_m(d, b, c\vert q) 
}
\eea

\noindent where $X_m(a, b, c| q)$ are 1-dimensional configuration 
sums, defined by 
\begin{equation}
X_m(a, b, c\vert  q)  = 
\sum_{configuration} q^{ \Phi_{m}(l_1, \cdots, l_{m+2})} \, .
\end{equation}
The sum in $X_m(a, b, c\vert  q)$ is taken over $l_2, \cdots, l_m$, with the
boundary heights fixed at $(l_1 = a, l_{m+1} = b, l_{m+2} = c)$. There are
$m + 2$ sites on a half-diagonal, counted starting from the centre
of the lattice, and proceeding towards an outer corner of the square
lattice.  

The 1-dimensional configuration sums
$X_m(a,b,c,\vert q)$  are precisely the objects  we are interested in.
To compute the local height probabilities in the thermodynamic limit,
one has to evaluate
\begin{equation}
          X(a, b, c \, \vert \, q )  =   
        \lim_{m\, even \rightarrow \infinity} X_m(a, b, c\, \vert\, q),  
\end{equation}
It is the configuration sum $X(a, b, c\vert q)$ 
that encodes the statistical mechanics of the model.

In \cite{ABF}, the 1-dimensional configuration sums were evaluated.
For example, in the case under consideration (regime III), it is
found that
$$
X_m(a,b,c\, \vert \, q)
=
q^{a(a-1)/4}\left( F_m(a,b,c) - F_m(-a,b,c) \right)
$$
where
$$
F_m(a,b,c)=
\sum_{n=-\infty}^{\infty}
q^{n(L-1)(nL-a) +[bc-(2nL-a)(b+c-1)]/4}
{m \brack {1\over 2}(m+a-b) - nL}_q
$$
In the thermodynamic limit, one obtains
\begin{equation}
X(a,b,c|q) 
=
{1\over \varphi(q)}
q^{bc/4}
\Delta(a,{1\over 2}(b+c-1)\, ;\, q)
\end{equation}
where
$$
\varphi(q)=\prod_{n\ge 1}(1-q^n)
$$
and
$$
\Delta(a,d\, ;\, q)
=
\sum_{n=-\infty}^\infty
q^{L(L-1)n^2+Ldn +a(a-1)/4}
\left(
q^{-(L-1)an-ad/2} - q^{(L-1)an+ad/2}
\right)
\, .
$$

In the case of the ABF models, the 1-dimensional configurations
$q$-counted by $X_m(a, b, c| q)$ can be represented graphically 
in terms of lattice paths (or Dyck graphs), defined
by the points of coordinates $(k-1,l_k-1)$ $(k=1,\ldots,m)$ in
the plane: 
\bigskip

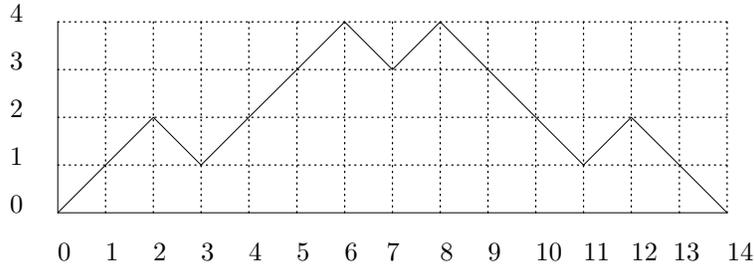
\begin{center}
\unitlength=0.00625in
\begin{picture}(612,224)(0,-10)
\path(40,200)(40,40)
\path(40,40)(80,80)(120,120)
	(160,80)(200,120)(240,160)
	(280,200)(320,160)(360,200)
	(400,160)(440,120)(480,80)
	(520,120)(560,80)(600,40)
\path(40,40)(600,40)
\dottedline{5}(40,80)(600,80)
\dottedline{5}(40,120)(600,120)
\dottedline{5}(40,160)(600,160)
\dottedline{5}(40,200)(600,200)
\dottedline{5}(600,200)(600,40)
\dottedline{5}(80,200)(80,40)
\dottedline{5}(120,200)(120,40)
\dottedline{5}(160,200)(160,40)
\dottedline{5}(200,200)(200,40)
\dottedline{5}(240,200)(240,40)
\dottedline{5}(280,200)(280,40)
\dottedline{5}(320,200)(320,40)
\dottedline{5}(360,200)(360,40)
\dottedline{5}(400,200)(400,40)
\dottedline{5}(440,200)(440,40)
\dottedline{5}(480,200)(480,40)
\dottedline{5}(520,200)(520,40)
\dottedline{5}(560,200)(560,40)
\put(0,40){\makebox(0,0)[lb]{\raisebox{0pt}[0pt][0pt]{\shortstack[l]{{\tenrm 0}}}}}
\put(0,80){\makebox(0,0)[lb]{\raisebox{0pt}[0pt][0pt]{\shortstack[l]{{\tenrm 1}}}}}
\put(0,120){\makebox(0,0)[lb]{\raisebox{0pt}[0pt][0pt]{\shortstack[l]{{\tenrm 2}}}}}
\put(0,160){\makebox(0,0)[lb]{\raisebox{0pt}[0pt][0pt]{\shortstack[l]{{\tenrm 3}}}}}
\put(0,200){\makebox(0,0)[lb]{\raisebox{0pt}[0pt][0pt]{\shortstack[l]{{\tenrm 4}}}}}
\put(40,0){\makebox(0,0)[lb]{\raisebox{0pt}[0pt][0pt]{\shortstack[l]{{\tenrm 0}}}}}
\put(80,0){\makebox(0,0)[lb]{\raisebox{0pt}[0pt][0pt]{\shortstack[l]{{\tenrm 1}}}}}
\put(120,0){\makebox(0,0)[lb]{\raisebox{0pt}[0pt][0pt]{\shortstack[l]{{\tenrm 2}}}}}
\put(160,0){\makebox(0,0)[lb]{\raisebox{0pt}[0pt][0pt]{\shortstack[l]{{\tenrm 3}}}}}
\put(200,0){\makebox(0,0)[lb]{\raisebox{0pt}[0pt][0pt]{\shortstack[l]{{\tenrm 4}}}}}
\put(240,0){\makebox(0,0)[lb]{\raisebox{0pt}[0pt][0pt]{\shortstack[l]{{\tenrm 5}}}}}
\put(280,0){\makebox(0,0)[lb]{\raisebox{0pt}[0pt][0pt]{\shortstack[l]{{\tenrm 6}}}}}
\put(315,0){\makebox(0,0)[lb]{\raisebox{0pt}[0pt][0pt]{\shortstack[l]{{\tenrm 7}}}}}
\put(360,0){\makebox(0,0)[lb]{\raisebox{0pt}[0pt][0pt]{\shortstack[l]{{\tenrm 8}}}}}
\put(400,0){\makebox(0,0)[lb]{\raisebox{0pt}[0pt][0pt]{\shortstack[l]{{\tenrm 9}}}}}
\put(440,0){\makebox(0,0)[lb]{\raisebox{0pt}[0pt][0pt]{\shortstack[l]{{\tenrm 10}}}}}
\put(480,0){\makebox(0,0)[lb]{\raisebox{0pt}[0pt][0pt]{\shortstack[l]{{\tenrm 11}}}}}
\put(520,0){\makebox(0,0)[lb]{\raisebox{0pt}[0pt][0pt]{\shortstack[l]{{\tenrm 12}}}}}
\put(555,0){\makebox(0,0)[lb]{\raisebox{0pt}[0pt][0pt]{\shortstack[l]{{\tenrm 13}}}}}
\put(600,0){\makebox(0,0)[lb]{\raisebox{0pt}[0pt][0pt]{\shortstack[l]{{\tenrm 14}}}}}
\end{picture}
\end{center}
\subsection{From 1-dimensional configurations to Virasoro
characters}  

The 1-dimensional configurations of the ABF models can be  
regarded as combinatorial objects with shape-dependent weights.
On that basis alone ABF evaluated the 1-dimensional sums, and
obtained the result in terms of $q$-series.

The Kyoto group made the fundamental observation that the generating 
functions of the ABF 1-dimensional sums are the characters of highest 
weight modules of Virasoro algebras \cite{DJMO1}. 
To be more precise, these are 
precisely the characters of the Virasoro modules arising in
the discrete unitary conformal field theories discovered  
by Belavin, Polyakov, and Zamolodchikov \cite{BPZ}.

\subsection{From Virasoro characters to branching functions}

The Kyoto group also observed that the ABF paths are intrinsically 
related to representation theory: they can be regarded as walks on 
a finite section of the weight lattice of the Lie algebra $\slchap_2$. 
This connection will be explained in more detail in the forthcoming section. 

Furthermore, they showed that the Virasoro characters corresponding
to configuration sums
can be obtained 
as the branching coefficients of products of affine characters,
as in the coset construction of Goddard-Kent-Olive
\cite{GKO}.

This will also be explained in more details in latter sections.
What we wish to point out here is that, once the connection with
$\slchap_2$ is made, it is natural to look for models that are related 
to higher rank Lie algebras. These models do exist \cite{JMO}.


\section{Combinatorics of solvable lattice models}

The $\slchap_n$ analogues of the ABF models were obtained in \cite{JMO}. 
For lack of space, we will not go into the definition of these models. 
Suffice it to say that they can be solved using the CTM, and that the 
result of the application of the CTM are configurations that can be 
interpreted as walks on restricted sections of the weight lattice of 
$\slchap_n$ as will be explained below.

To represent these generalised 1-dimensional configurations\footnote{The 
can still be regarded as one dimensional configurations, in the sense 
that the can be represented as a linear sequence of sets of integers.}   
the ABF-type paths are no longer adequate,
and one has to work in terms of `coloured' partitions.  

\subsection{Representation theory of $\slchap_n$}

\subsubsection{The affine Lie algebra $\slchap_n$}\label{Fundamentals}
Let $n\ge2$ and let $\h$ be a 
$(n+1)$-dimensional vector space over $\Q$ with basis 
$\{h_0,h_1,\ldots ,h_{n-1},D\}$.
We let 
$\{\Lambda_0,\Lambda_1,\ldots ,\Lambda_{n-1},\delta\}$
be the corresponding dual basis of $\h^*$, the dual space of $\h$.
That is,
$$
\<\Lambda_i,h_j\> = \delta_{ij}, \quad
\<\Lambda_i,D\> = 0,\quad
\<\delta, h_i\> = 0, \quad
\<\delta, D\> = 1 .
$$
It will be convenient to extend the index set so that
$\Lambda_{i}=\Lambda_{(i\mod n)}$ for all $i\in\Z$.
Then, for all $i\in\Z$, we set
$\epsilon_i=\Lambda_{i+1}-\Lambda_i$ and
$\alpha_i=2\Lambda_i-\Lambda_{i-1}-\Lambda_{i+1} + \delta_{i0} \,\delta$,
where $\delta_{ij}=1$ if $i-j\mod n=0$ and $\delta_{ij}=0$ otherwise.

The $n\times n$ matrix $[\<\alpha_i,h_j\>]$ is called 
the generalised Cartan matrix of type $A_{n-1}^{(1)}$.
The associated (affine) Kac-Moody algebra is denoted by $\slchap_n$.
It is defined as the algebra generated by
$D$ and $e_i$, $f_i$, $h_i$ for $0\le i<n$, subject to the relations:
$$
[h_i,h_j]=0;\quad [h_i,D]=0;
$$ $$
[h_i,e_j]=\<\alpha_j,h_i\>e_j; \quad [D,e_j]=\delta_{j0}e_j;
$$ $$
[h_i,f_j]=-\<\alpha_j,h_i\>f_j; \quad [D,f_j]=-\delta_{j0}f_j;
$$ $$
[e_i,f_j]=\delta_{ij}h_i;
$$ $$
({\rm ad}\,e_i)^{1-\<\alpha_j,h_i\>} e_j=0 \quad (i\ne j);
$$ $$
({\rm ad}\,f_i)^{1-\<\alpha_j,h_i\>} f_j=0 \quad (i\ne j),
$$
where $({\rm ad}\,a)\,b=[a,b]$.

The $\Lambda_i$ are known as the fundamental weights, the
$\epsilon_i$ the fundamental vectors, and the $\alpha_i$ the simple roots
of $\slchap_n$.

\subsubsection{Weight lattices and root lattices}\label{WeightLattices}

The root lattice of $\slchap_n$ is 
$$
Q=\bigoplus_{i=0}^{n-1} \Z\alpha_i.
$$
Its weight lattice  is 
$$
P=\Z\delta \oplus \bigoplus_{i=0}^{n-1} \Z\Lambda_i.
$$
Elements of $P$ are known as weights.
The dual weight lattice of $\slchap_n$ is defined to be
$$
P^\vee = \Z\,D\oplus \bigoplus_{i=0}^{n-1} \Z\,h_i,
$$
For $l\in\Z$, $P_l$ is defined by
$$
P_l=\left\{\sum_{i=0}^{n-1}a_i\Lambda_i\mid a_i\in\Z \hbox{ for }
i=0,1,\ldots,n-1; \sum_{i=0}^{n-1}a_i=l\right\}.
$$
For $l\in\N$, $P^+_l$ is defined by
$$
P_l^+=\left\{\sum_{i=0}^{n-1}a_i\Lambda_i\in P_l\mid a_i\ge0\right\}.
$$
Then define $P^+=\bigcup_{l\ge0} P^+_l$.

Strictly speaking, $P_l$ is the level $l$ integral weight
lattice of $\slchap_n{}'$, which is the subalgebra of
$\slchap_n$ generated by $e_i,h_i,f_i$ (without $D$).

\subsubsection{$\slchap_n$-modules} \label{slnModules}

Let $M$ be an $\slchap_n$-module.
For each weight $\Lambda \in P$, the subspace of $M$ defined by
$$
M_\Lambda = \{ v\in M \ | \ h\,v = \<\Lambda,h\> \,v, \ h\in\h \}
$$
is called the weight space of weight $\Lambda$ of $M$.
If $v\in M_\Lambda$ then we write $\wt(v)=\Lambda$, and
call such $v$ a weight vector of weight $\Lambda$.

The module $M$ is said to be integrable if:
\renewcommand{\labelenumii}{\roman{enumii}.}
\begin{enumerate}
\item $M=\bigoplus_{\Lambda \in P} M_\Lambda$;
\item ${\rm dim}\,M_\Lambda < \infty$ for each $\Lambda \in P$;
\item for each $i= 0,1,\ldots ,n-1,$ $M$ decomposes into a direct sum
of finite dimensional ${\goth g}_i$-modules, where ${\goth g}_i$ denotes 
the subalgebra of $\slchap_n$ generated by $e_i$, $f_i$, and $h_i$.
\end{enumerate}

\noindent
A weight vector $v\in M$ for which $e_iv=0$ for all
$i=0,1,\ldots,n-1$, is said to be a highest weight vector.
If there exists a highest weight vector $v\in M$ such that
$M = U(\slchap_n)\, v$, then $M$ is said to be a highest weight module.
Then the weight of $v$ is called the highest weight of $M$.

\noindent
If $l\ge0$ then for all $\Lambda\in P^+_l$, there exists
(up to equivalence)
a unique irreducible integrable
$\slchap_n$-module $V(\Lambda)$ with highest weight $\Lambda$
\cite{KacBook}.

\subsubsection{The Virasoro algebra}

The Virasoro algebra is the complex Lie algebra $Vir$ generated
by elements $L_n$ ($n\in {\bf Z}$) and $c$ subject to the
relations
$$
[L_m\, ,\, L_n] = (m-n)L_{m+n} +{1\over 12}(m^3-m)\, c \, \delta_{m+n, 0}
$$
and 
$$
[c\, ,\, L_n] = 0 \ .
$$
It is the universal central extension of the Witt algebra, {\it i.e. }the 
Lie algebra generated by the differential operators $L_n=-z^{n+1}{d\over dz}$.  
On an irreducible representation, $c$ is a scalar operator, whose eigenvalue 
$c$ is called {\sl the central charge}. The eigenvalue of $L_0$ on the highest 
weight vector is called the conformal dimension $h$.

By using the Sugawara construction
(see \cite{KacBook}), one can construct an action
of $Vir$ on any representation of $\slchap_n$ of fixed level $l$.
This can be done in particular for the tensor product of
two irreducible representations \cite{GKO}.

\subsubsection{Branching functions}

Let $\Lambda'$, $\Lambda''$ be dominant integral weights of
respective levels $l'$, $l''$.
The tensor product of $V(\Lambda')$ by $V(\Lambda'')$ can be
written as the finite direct sum
\begin{equation}
V(\Lambda')\otimes V(\Lambda'')
=
\bigoplus_{\Lambda\in P^+_l} \Omega^\Lambda_{\Lambda'\Lambda''}
\otimes
V(\Lambda)
\end{equation}
where $l=l'+l''$ and
$\Omega^\Lambda_{\Lambda'\Lambda''}$ is the space of highest
weight vectors whose weights are congruent to $\Lambda$ modulo
$\delta$. This is the decomposition into isotypic components
for the action of $\slchap_n{}'$ ({\it i.e.}  $\slchap_n$ without
the generator $D$).

The action of $Vir$ on $V(\Lambda')\otimes V(\Lambda'')$
defined by the
coset construction of \cite{GKO}
commutes with the action of $\slchap_n{}'$,
so that each $\Omega^\Lambda_{\Lambda'\Lambda''}$
is a representation of $Vir$.

The formal character
$$
b^\Lambda_{\Lambda'\Lambda''}(q)
= {\rm ch\,}\Omega^\Lambda_{\Lambda'\Lambda''}
=
\sum_k  q^k\dim (\Omega^\Lambda_{\Lambda'\Lambda''})_{\Lambda-k\delta}
$$
is called a branching function (this definition differs from the
usual Virasoro character by a factor $q^h$).

For $n=2$, these Virasoro characters are precisely those of the 
discrete unitary series of irreducible representations. For $n>2$, 
the spaces $\Omega^\Lambda_{\Lambda'\Lambda''}$ are generally not 
irreducible under the Virasoro algebra. They are irreducible 
representations of the ${\cal W}_n$-algebras, which are less 
well-understood.

The ABF models are related to the case $n=2$. Let $V_{j,l}$ be the 
irreducible level $l$ representation $V((l-j)\Lambda_0+j\Lambda_1)$ 
of $\slchap_2$, and let $b_{j_1j_2j_3}(q)$ be the character of 
$\Omega_{j_1j_2j_3}$, where
$$
V_{j_1,L-3}\otimes V_{j_2,1}
=
\bigoplus_{j_3} \Omega_{j_1j_2j_3}\otimes V_{j_3,L-2}
$$
Then, the 1-dimensional configuration sums are given by
$$
X(a,b,c|q)
=q^\gamma b_{j_1j_2j_3}(q)
$$
where 
$$
j_1={1\over 2}(b+c-1)-1,\quad
j_2 = {1\over 2}(b-c+1), \quad
j_3 = a-1
$$
and $\gamma$ is a certain rational number \cite{DJMO1,DJMO2}.

Let $L(c,h)$ be the irreducible representation of the
Virasoro algebra with central charge $c$
and conformal weight $h$.
In the $\slchap_2$ case, the decomposition of the tensor product reads
$$
V_{j,l}\otimes V_{\varepsilon,1}
=
\bigoplus_s
L(c,h_{rs})\otimes V_{s-1,l+1}
$$
where $r=j+1$ and the sum runs over all $s$ such that
$1\le s \le l+2$, $r-s$ even if $\varepsilon=0$ or
$r-s$ odd if $\varepsilon=1$. The central charge
$c$ is given by
$$
c= 1 - {6\over (l+2)(l+3)}
$$
and 
$$
h_{rs}={ [(l+3)r-(l+2)s]^2-1 \over 4(l+2)(l+3)}
$$
The character $\chi_{r,s}(q)$ of $L(c,h_{rs})$ is given by the
Rocha-Caridi formula \cite{RC}
$$
\chi_{r,s}(q)
=
{1\over \varphi(q)}
\sum_{k=-\infty}^\infty
\left(
-q^{a(k)}+q^{b(k)}
\right)
$$
where, setting $m=l+2$
\bea
a(k)& = &
{ [2m(m+1)k+(m+1)r+ms]^2-1 \over 4m(m+1)}\\ \nonumber
b(k) &=&
{ [2m(m+1)k +(m+1)r - ms]^2 -1 \over 4m(m+1) } \ .\\ \nonumber
\eea

\subsection{Paths}

As explained above, the notion of a {\em path} emerged from 
the corner transfer matrix method. Roughly speaking, a path 
corresponds to an eigenvector of the corner transfer matrix 
at the absolute temperature $q=0$. It is defined as a certain 
sequence of points on the weight lattice of $\slchap_n$.

\begin{definition}[Paths]\label{PathDef}
\hfil
\begin{enumerate}
\item A level $l$ path is a sequence
$$
p=(p_0,p_1,p_2,\ldots),
$$
for which $p_k\in P_l$ for all $k\ge0$.
The set of all level $l$ paths is denoted ${\cal P}_l$.
\item If $p\in{\cal P}_{l}$ and $\mu\in P_m$,
define the level $l-m$ path $p-\mu$ by
$$
p-\mu=(p_0-\mu,p_1-\mu,p_2-\mu,\ldots).
$$
\end{enumerate}
\end{definition}

\begin{definition}[Unrestricted paths]\label{vPathDef}
\hfil
\begin{enumerate}
\item A $\Lambda_0$-path $p$ is a level 1 path for which,
for all $k\ge0$, there exists $\gamma(k)$ satisfying
$0\le\gamma(k)<n$, such that $p_{k+1}-p_k=\epsilon_{\gamma(k)}$,
and for which for some $k_*\ge0$, $p_k=\Lambda_{k}$
for all $k\ge k_*$.
\item The set of all $\Lambda_0$-paths is denoted ${\cal P}(\Lambda_0)$.
\item The length of $p\in{\cal P}(\Lambda_0)$ is the minimal $k_*$
for which $p_k=\Lambda_{k}$ for all $k\ge k_*$.
We then define $l(p)=k_*$.
\item The ground state of ${\cal P}(\Lambda_0)$ is the unique
path $\overline p\in{\cal P}(\Lambda_0)$ having length 0.
Thus $\overline p_k=\Lambda_{k}$ for all $k\ge0$.
\end{enumerate}
\end{definition}

\noindent
Note that a $\Lambda_0$-path $p$ is completely determined by the sequence
$\gamma(0)$, $\gamma(1)$, $\gamma(2),\ldots$.
The adjective {\it unrestricted} will be used when appropriate
to distinguish the paths defined here from restricted paths
which are defined later.

\begin{example}\label{PathEx}{\rm
Let $n=2$.
The following diagram depicts the $\Lambda_0$-path
$(\Lambda_0$, $\Lambda_1$, $2\Lambda_1-\Lambda_0$,
$3\Lambda_1-2\Lambda_0,2\Lambda_1-\Lambda_0,
\Lambda_1,2\Lambda_1-\Lambda_0,
\Lambda_1,\Lambda_0,
2\Lambda_0-\Lambda_1,\Lambda_0,
\Lambda_1,\Lambda_0$,
$\Lambda_1$, $\Lambda_0,\ldots)$.

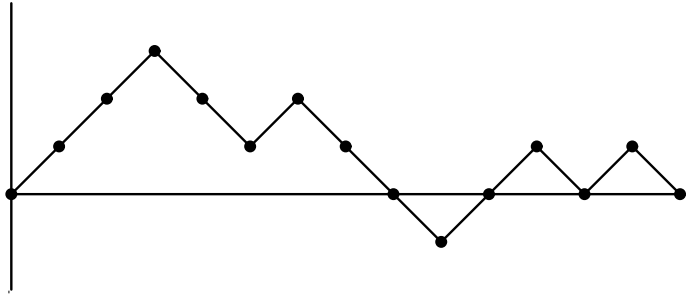
\medskip
\centerline{\epsfbox{path.eps}}
\medskip

\noindent This path has length 10.
}
\end{example}

With each $\Lambda_0$-path, we can associate an {\it energy}.

\begin{definition}\label{EnergyDef}
Let $p\in{\cal P}(\Lambda_0)$,
and $\gamma(k)$ be defined by \ref{vPathDef}.1.
Let $\overline\gamma(k)$ be the corresponding values for the
ground state $\overline p\in{\cal P}(\Lambda_0)$:
$\overline\gamma(k)=k\bmod n$ for all $k\ge0$.
\begin{enumerate}
\item The energy of $p$ is defined to be
$$
E(p)=\sum_{k=1}^\infty
k\Big(H\big(\gamma(k-1),\gamma(k)\big)
-H\big(\overline\gamma(k-1),\overline\gamma(k)\big)\Big),
$$
where, for $0\le a,b<n$, we define
$$
H(a,b)=\cases{0&if $a<b$;\cr 1&if $a\ge b$.\cr}
$$
\item The weight of $p$ is then defined to be
$$
\wt(p)=p_0-E(p)\delta,
$$
where $\delta$ is the null root of $\slchap_n$.
\end{enumerate}
\end{definition}

\noindent
Note that the energy $E(p)$ of each path $p\in{\cal P}(\Lambda_0)$
is finite because $\gamma(k)=\overline\gamma(k)$ for all $k\ge l(p)$.

\subsection{Partitions and paths} \label{subsec:n-reg}

A partition $\lambda$ of $m$ is a sequence
$(\lambda_1,\lambda_2,\ldots,\lambda_r)$ of positive integers
such that $\lambda_1\ge\lambda_2\ge\cdots\ge\lambda_r$
and $\sum_{i=1}^r \lambda_i=m$. For $i>r$, we define $\lambda_i=0$.
Occasionally the notation
$\lambda=(\lambda_1^{a_1},\lambda_2^{a_2},\ldots,\lambda_r^{a_r})$
will be used to denote a partition, where for $i=1,2,\ldots,r$,
the part $\lambda_i$ occurs in $\lambda$ with a multiplicity $a_i$.

The (finite) set of all partitions of $m$ is
denoted $\Pi(m)$. Then define $\Pi=\cup_{m\ge0}\Pi(m)$.

The diagram $F^\lambda$ associated with $\lambda\in\Pi(m)$ consists 
of $m$ nodes (or boxes) arranged in $r$ left adjusted rows. The number 
of nodes in the $i$th row is $\lambda_i$. For example:
$$
F^{(4,3,1)}=\smyoungd{
\multispan9\hrulefill\cr &&&&&&&&\cr
\multispan9\hrulefill\cr &&&&&&\cr
\multispan7\hrulefill\cr &&\cr
\multispan3\hrulefill\cr}.
$$
The partition $\lambda^\prime$ conjugate to the
partition $\lambda\in\Pi(m)$ is obtained by setting
$\lambda^\prime_i$ to be the length of the $i$th column
of $F^\lambda$ for $1\le i\le\lambda_1$ (reading from the left).

A partition $\lambda\in\Pi(m)$ is said to be $n$-regular
if no part appears $n$ or more times.
Equivalently: $\lambda^\prime_i-\lambda^\prime_{i+1}<n$ for
all $i\ge1$.

The set of $n$-regular partitions of $m$ is
denoted $\Pi_n(m)$. Then define $\Pi_n=\cup_{m\ge0}\Pi_n(m)$.

\noindent
In the following definition, we associate an $n$-regular partition with
each path $p\in{\cal P}(\Lambda_0)$.

\begin{definition}[Highest-lift of a path]
Let $p\in{\cal P}(\Lambda_0)$ and $k_*=l(p)$.
For $k\ge k_*$, let $t_k=0$.
Then for $k=k_*-1,k_*-2,\ldots,0$, recursively calculate $t_k$
from $t_{k+1}$ by setting $t_k=k-\gamma(k)\pmod n$ with
$0\le t_k-t_{k+1}<n$.
The partition $\lambda=\lambda(p)$ is then defined such that
its conjugate is
$$
\lambda^\prime=(t_0,t_1,\ldots,t_{k_*-1}).
$$
The partition $\lambda(p)$ is known as the highest-lift of $p$.
\end{definition}

\begin{example}\label{PathtoPartitionEx}{\rm
Consider the path described in Example \ref{PathEx}.
Here we calculate that for $k=0,1,..,9,$ the values of
$\gamma(k)$ are $0,0,0,1,1,0,1,1,1,0,$ respectively.
Hence, modulo 2, the respective values of $t_k$ are
$0,1,0,0,1,1,1,0,1,1.$
Since the path has length 10, we set $t_{10}=0$.
Then, using $0\le t_k-t_{k+1}<n$, we obtain the values
$t_9=1$, $t_8=1$, $t_7=2$, $t_6=3$, $t_5=3$, $t_4=3$, $t_3=4$,
$t_2=4$, $t_1=5$, $t_0=6$.
This yields the following partition: 
$$
\smyoungd{
\multispan{21}\hrulefill\cr &&&&&&&&&&&&&&&&&&&&\cr
\multispan{21}\hrulefill\cr &&&&&&&&&&&&&&&&\cr
\multispan{17}\hrulefill\cr &&&&&&&&&&&&&&\cr
\multispan{15}\hrulefill\cr &&&&&&&&\cr
\multispan{9}\hrulefill\cr &&&&\cr
\multispan{5}\hrulefill\cr &&\cr
\multispan{3}\hrulefill\cr}\;,
$$
and hence the highest-lift partition $(10,8,7,4,2,1)$.
}
\end{example}

\noindent
Note that $\lambda(\overline p)$ is the empty partition
and that each $\lambda(p)$ is $n$-regular.
Moreover, each $n$-regular partition corresponds to a unique
element of ${\cal P}(\Lambda_0)$.
Thus the sets ${\cal P}(\Lambda_0)$ and $\Pi_n$ are in bijection.

It is often useful to fill each node of the partition
corresponding to a path $p$ with an integer.

\begin{definition}[Colouring]\label{ColouringDef}
A partition $\lambda$ is said to be coloured with $v$ if,
for $i=1,2,\ldots$, and $j=1,2,\ldots,\lambda_i$, the
$j$th node in the $i$th row is filled with the value
$j-i+v\mod n$. The value occupying a node is known as
its colour charge.
The colour charge $v$ of the leftmost node in the first row
is referred to as the colour of the partition.
The weight of the partition is $\Lambda_v$.
\end{definition}

\noindent
For example, on colouring the partition $\lambda=(5,5,4,1,1)$
with $0$ in the case $n=3$, we obtain:
$$
\smyoungd{
\multispan{11}\hrulefill\cr &0&&1&&2&&0&&1&\cr
\multispan{11}\hrulefill\cr &2&&0&&1&&2&&0&\cr
\multispan{11}\hrulefill\cr &1&&2&&0&&1&\cr
\multispan{9}\hrulefill\cr &0&\cr
\multispan{3}\hrulefill\cr &2&\cr
\multispan{3}\hrulefill\cr }.
$$

\begin{definition}\label{EwtDef}
Colour the partition $\lambda$ with 0 and, for $0\le i<n$, let
$m_i$ be the multiplicity of the colour $i$ in $\lambda$.
Then define:
\begin{enumerate}
\item $E(\lambda)=m_0$;
\item $\wt(\lambda)=\La_0-\sum_{i=0}^{n-1}m_i\alpha_i$.
\end{enumerate}
\end{definition}

\noindent When $\lambda$ is interpreted as a basis vector
$s_\lambda$ of the Fock space (see Section \ref{Sec4.2})
$\wt(\lambda)$ is the $\Sl_n$-weight of $s_\lambda$,
and $E(\lambda)$ its degree in the homogeneous gradation.
It is straightforward to show that the energy and weight
of a path coincide with those of the corresponding partition:

\begin{lemma}\label{EwtLemma}
Let $p\in{\cal P}(\Lambda_0)$ and $\lambda=\lambda(p)$, the corresponding
partition.
\begin{enumerate}
\item $E(p)=E(\lambda)$;
\item $\wt(p)=\wt(\lambda)$.
\end{enumerate}
\end{lemma}

The next result, which was first given in \cite{DJKMOl}, states that
the character of the basic representation $V(\Lambda_0)$ of $\slchap_n$,
amounts to enumerating the $n$-regular partitions, taking
account of their colouring.

\begin{theorem}\label{FormalChar} {\rm \cite{DJKMOl,JMMO}}
The formal character of the basic representation $V(\Lambda_0)$ of
$\slchap_n$ is
$$
\ch V(\Lambda_0)=\sum_{\lambda\in\Pi_n}
e^{\wt(\lambda)}.
$$
\end{theorem}

\noindent
The principally specialised character of $V(\Lambda)$, defined
by
\begin{equation}\label{PrincipalCharDef}
{\rm Pr\,}\ch V(\Lambda)
=\sum_{m_0,\ldots,m_{n-1}}
\dim V(\Lambda)_{\Lambda-m_0\alpha_0-\cdots-m_{n-1}\alpha_{n-1}}
q^{m_0+\cdots+m_{n-1}},
\end{equation}
then follows in the $\Lambda=\Lambda_0$ case, as a direct corollary:

\begin{theorem}\label{PrincipalChar} {\rm\cite{KacBook}}
$$
{\rm Pr\,}\ch V(\Lambda_0)
=\prod_{\scriptstyle j\ge1\atop\scriptstyle j\not\equiv0\,
\mbox{\scriptsize mod}\,n}
\frac1{1-q^j}.
$$
\end{theorem}
\Proof By Theorem \ref{FormalChar} and Definition \ref{EwtDef}.2,
the coefficient of $q^N$ on the right side of
(\ref{PrincipalCharDef}) is given by the number of $n$-regular 
partitions of $N$. The generating function for this number 
\cite{JK} then gives the result.

\medskip\noindent
We will later describe a realisation of $V(\Lambda_0)$ that has a
basis naturally indexed by the set of $n$-regular partitions $\Pi_n$.

\subsection{Restricted paths}

In Definition \ref{vPathDef}, we defined a set of unrestricted paths
${\cal P}(\Lambda_0)$. We now use these to define restricted paths.

\begin{definition}[Restricted paths]\label{rPathDef}
Let $\mu$ be a dominant integral weight of level $l$.

\begin{enumerate}
\item A path $p$ is said to be $(\mu,\Lambda_0)$-restricted
      if $p-\mu\in{\cal P}(\Lambda_0)$ and $p_k\in P^+_{l+1}$ 
      for all $k\ge0$.  
      
\item The set of $(\mu,\Lambda_0)$-restricted paths is denoted
      ${\cal P}(\mu,\Lambda_0)$.  
      
\item The length $l(p)$ of $p\in{\cal P}(\mu,\Lambda_0)$ is
      defined to be $l(p-\mu)$.

\item The ground state $\overline p\in{\cal P}(\mu,\Lambda_0)$
      is defined to be $\overline p=\overline p^\prime+\mu$
      where $\overline p^\prime$ is the ground state of 
      ${\cal P}(\Lambda_0)$.

Thus here, $\overline p_k=\mu+\Lambda_k$ for all $k\ge0$.

\end{enumerate}
\end{definition}

\noindent
The elements of ${\cal P}(\mu,\Lambda_0)$ are called {\em restricted 
paths}, since they describe 'walks' on a restricted segment of a
weight lattice. They also label the 1-dimensional configurations of 
certain restricted solid-on-solid models. Since
$\{p - \mu | p \in {\cal P}$ $(\mu,\Lambda_0)\}$ $\subset$ ${\cal P}$ 
$(\Lambda_0)$, a partition may readily be associated with each 
restricted path.

\begin{definition}[Restricted partitions]
\hfil
\begin{enumerate}

\item Let $p\in{\cal P}(\mu,\Lambda_0)$. Then the partition $\lambda(p)$ 
      associated with $p$ is defined to be $\lambda(p-\mu)$.

\item Define 
${\cal Y}(\mu,\Lambda_0)
=\left\{\lambda(p)\mid p\in{\cal P}(\mu,\Lambda_0)\right\}$.

\end{enumerate}
\end{definition}

\noindent
Clearly ${\cal Y}(\mu,\Lambda_0)\subset\Pi_n$. In the case where $\mu$ is 
of level one, this set of partitions was first characterised in \cite{FOW}.

\subsection{The FOW conditions}

\begin{theorem}[FOW theorem] \label{FOWTheorem} {\rm\cite{FOW}}
Let 
$\lambda =(\lambda_1^{a_1}\!,\lambda_2^{a_2}\!,\ldots,\lambda_r^{a_r})
\!\in\Pi_n$,
where 
$\lambda_1>\lambda_2>\cdots>\lambda_r>0$ 
and
$0<a_i<n$ for $i=1,2,\ldots,r$.
Then 
$\lambda\in{\cal Y}(\Lambda_j,\Lambda_0)$ 
if and only if either $r=0$, or
$$
a_i+\lambda_i-\lambda_{i+1}+a_{i+1}=0 \pmod n,
$$
for $i=1,2,\ldots,r-1$,  and $j=\lambda_1-a_1\,(\mod n)$.
\end{theorem}

Note that the value of $\lambda_i-\lambda_{i+1}$ is the length of
a horizontal edge of the corresponding partition (it is the $i$th 
horizontal edge from the top, not counting the top border of the 
diagram). $a_i$ and $a_{i+1}$ are the lengths of its neighbouring 
vertical edges. Thus, we must check whether the sum of each horizontal 
edge and its neighbouring vertical edges is a multiple of $n$. This 
makes identifying these partitions particularly straightforward.

\begin{example}\label{FOWEx}{\rm
Consider the following partition, which is coloured with $0$
in the case $n=3$:
$$\lambda\:=\:
\smyoungd{
\multispan{27}\hrulefill\cr 
&0&&1&&2&&0&&1&&2&&0&&1&&2&&0&&1&&2&&0&\cr
\multispan{27}\hrulefill\cr
&2&&0&&1&&2&&0&&1&&2&&0&&1&&2&&0&&1&&2&\cr
\multispan{27}\hrulefill\cr
&1&&2&&0&&1&&2&&0&&1&&2&&0&&1&\cr
\multispan{21}\hrulefill\cr
&0&&1&&2&&0&&1&&2&\cr
\multispan{13}\hrulefill\cr
&2&&0&&1&&2&&0&\cr
\multispan{11}\hrulefill\cr
&1&&2&&0&&1&\cr
\multispan{9}\hrulefill\cr
&0&\cr
\multispan{3}\hrulefill\cr
&2&\cr
\multispan{3}\hrulefill\cr}\;.
$$
Starting from the top right of this partition, the lengths of
the vertical and horizontal edges taken in order along the border
are 2,3,1,4,1,1,1,1,1,3,2,1. The above description requires that
we sum the three successive values beginning at the 1st, 3rd, 5th,
etc. We obtain the values 6,6,3,3,6.
Since each is divisible by $3$, we conclude that
$\lambda\in{\cal Y}(\Lambda_j,\Lambda_0)$ for a certain $j$.
Theorem \ref{FOWTheorem} also gives
$j=(\lambda_1-a_1)\,\mod n=(13-2)\,\mod 3=2$.

The above vertical-horizontal-vertical edge lengths: 6,6,3,3,6;
are not all zero modulo $n\ne3$ and thus for $n\ne3$,
$\lambda\notin{\cal Y}(\Lambda_j,\Lambda_0)$ for all $j$.

As a further example, consider the partition $\mu=(10,8,7,4,2,1)$
which was produced in Example \ref{PathtoPartitionEx}.
In this case, we find the sequence 1,1,1, of vertical-horizontal-vertical
edge lengths. Since their sum is $3\ne0\,(\mod 2)$, we conclude
that in the case $n=2$,
$\mu\notin{\cal Y}(\Lambda_j,\Lambda_0)$ for all $j$.
Hence, the path of Example \ref{PathEx}
is not an element of ${\cal P}(\Lambda_j,\Lambda_0)$ for any $j$.
}
\end{example}

\noindent
It is worth noting that the empty partition $\emptyset$ satisfies
the conditions of Theorem \ref{FOWTheorem} trivially and hence
$\emptyset\in{\cal Y}(\Lambda_j,\Lambda_0)$ for all $j$ and all $n$.
Correspondingly, the appropriate ground state path
$\overline p\in{\cal P}(\Lambda_j,\Lambda_0)$.

The following result of \cite{JMMO} shows how the branching
function of the tensor product of two level one
$\slchap_n$ representations is obtained by the enumeration
of restricted partitions.           

\begin{theorem}{\rm \cite{JMMO}}\label{ThJMMO}
Let $0\le j<n$. Then for each weight $\Lambda$, the branching 
function $b^{\Lambda}_{\Lambda_j,\Lambda_0}(q)$ is given by:
$$
b^{\Lambda}_{\Lambda_j,\Lambda_0}
=\sum_{
\scriptstyle p\in{\cal P}(\Lambda_j,\Lambda_0)\atop
\scriptstyle p_0=\Lambda}
q^{E(p)}
=\sum_{
\scriptstyle\lambda \in{\cal Y}(\Lambda_j,\Lambda_0)\atop
\scriptstyle\wt(\lambda)=\Lambda -\Lambda_j\,(\mod\delta)}
q^{E(\lambda)}.
$$
\end{theorem}

\noindent
Here, the second equality results from the previously described
mapping between paths and partitions.

If we now define $FOW(n,j,k)$ to be the subset of
${\cal Y}(\Lambda_j,\Lambda_0)$ comprising those partitions $\lambda$
for which $\wt(\lambda)=\Lambda_k+\Lambda_{j-k}-\Lambda_j\,(\mod\delta)$,
we immediately obtain the following:

\begin{theorem}{\rm \cite{FOW}}
Let $0\le j<n$ and $0\le k\le(j-k)\,\mod n$.
Then the branching function
$b^{\Lambda_k+\Lambda_{j-k}}_{\Lambda_j,\Lambda_0}(q)$
is given by:
$$
b^{\Lambda_k+\Lambda_{j-k}}_{\Lambda_j,\Lambda_0}
=\sum_{\scriptstyle\lambda \in FOW(n,j,k)}
q^{E(\lambda)}.
$$
\end{theorem}

\subsection{Counting restricted paths}

For $L\ge0$, we set 

\begin{eqnarray}
 \P(\La_j, \La_0;L )&=& \{p \in {\cal P}(\La_j,\La_0) \mid l(p)\le L\}, \\
 b_{\La_j\,\La_0}^\La(q;L) &=& \sum_{p\in\P(\La_j,\La_0;L) 
\atop
p_0=\La}q^{E(p)},
\end{eqnarray}

\noindent
whereupon
$b_{\La_j\,\La_0}^{\La}(q;L)$ is a polynomial in $q$. 
In \cite{FOW}, a {\em constant sign q-series} for 
$b_{\La_j\,\La_0}^{\La}(q;L)$ was obtained. 
Since the expression in \cite{FOW} is in the principal picture, we 
need to adjust the overall power of $q$. We set $\La=\La_s+\La_t$ 
($0\le s\le t<n$). Let $C$ be the Cartan matrix of $\goth{sl}_n$, 
and $e_i$ be the $(n-1)$-dimensional unit vector 
$(0,\cd,0,\stackrel{i}{1},0\cd,0)^t$. We set $e_n=0$.
The expression reads as follows.

\begin{eqnarray}
b_{\La_j\,\La_0}^{\La}(q;L)&=&\sum_mq^{m^tC^{-1}m-m^tC^{-1}e_{s-t+n}+st/n}
\prod_{i=1}^{n-1}{l_i+m_i\brack m_i},\\
l&=&C^{-1}(Le_{n-1}+e_r+e_{s-t+n}-2m),
\end{eqnarray}

\noindent
where the sum is taken over all $m\in(\Zn)^{\times(n-1)}$ satisfying
$t+\sum_{i=1}^{n-1}im_i = 0\ (\mod n)$, and $r$ is determined from 
$L-(s+t)\equiv r$, $0<r\le n$.
Here, the notation means\footnote{The 
definition of $q$-binomial coefficient here is different from the one used
in \S 4.} 
$$
(q)_k=(1-q)(1-q^2)\cd(1-q^k),\quad{m\brack k}=\frac{(q)_m}{(q)_{m-k}(q)_k}.
$$
We refer to \cite{FOW} for details of the proof.
Taking the limit $L\rightarrow\infty$, we have
\[
b_{\La_j\,\La_0}^{\La}(q)=\sum_m\frac{q^{m^tC^{-1}m-m^tC^{-1}e_{s-t+n}+st/n}}
{\prod_{i=1}^{n-1}(q)_{m_i}},
\]
with the same restriction on $m$.


\goodbreak
\section{Modular representations of symmetric groups and Hecke algebras}\label{sect:2}

\subsection{The symmetric groups}\label{SGsection}

The symmetric group $\SG_m$ may be defined as the group
generated by the elements $s_i$, for $i=1,2,\ldots,m-1$,
subject to the relations:
\begin{eqnarray*}
&&s_is_{i+1}s_i=s_{i+1}s_is_{i+1};\\
&&s_is_j=s_js_i \qquad \vert{i-j}\vert>1;\\
&&s_i^2=1.
\end{eqnarray*}

\noindent
In the realisation of $\SG_m$ as the permutation group on the
set $\{1,2,\ldots,m\}$, the generator $s_i$ is identified with
the simple transposition $(i,i+1)$.
The subgroup generated by $s_i\ (i=1,\ldots m-2)$ is 
isomorphic to $\SG_{m-1}$, and we shall be concerned with
the problem of restricting representations of $\SG_m$
to this subgroup.
Over fields of characteristic 0, the representation theory of
$\SG_m$ has been understood since the beginning of the century.

\begin{theorem} Over fields of characteristic 0, the
inequivalent irreducible representations of $\SG_m$
are indexed by $\Pi(m)$, the set of partitions of $m$.
\end{theorem}

\noindent
There exist a number of ways to  calculate the
matrices of the irreducible representation of $\SG_m$
corresponding to a particular $\lambda\in\Pi(m)$.
One way is by means of the Specht module $S^\lambda$
(see \cite{JK} for example), whereby the resulting matrices
contain only integers.
On restricting to $\SG_{m-1}$, the
module $S^\lambda$ is no longer irreducible in general.
One has the following well-known branching rule.

\begin{theorem}\label{SpechtReduction}
Let $\lambda\in\Pi(m)$. Then
$$
S^\lambda \downarrow^{\SG_m}_{\SG_{m-1}}
=\bigoplus_{\mu\in{\cal R}(\lambda)} S^\mu,
$$
where ${\cal R}(\lambda)$ is the set of all partitions of $m-1$
that can be obtained by removing a single node from $\lambda$.
\end{theorem}

\noindent
The following lattice of partitions, which is known as Young's lattice,
results from joining a partition $\lambda\in\Pi(m)$ to a partition
$\mu\in\Pi(m-1)$, if $S^\mu$ occurs in
$S^\lambda \downarrow^{\SG_m}_{\SG_{m-1}}$:

\medskip
\centerline{\epsfbox{young.eps}}
\medskip
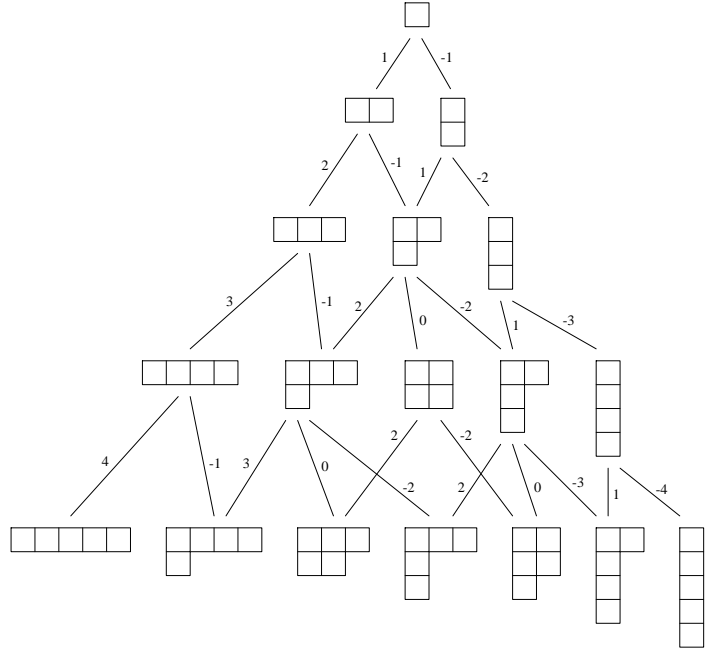
\centerline{Figure 1.}
\medskip

\noindent As indicated in Fig.~1, each edge is naturally labelled 
by an integer. 
If the edge corresponds to the removal of a node sitting on row $i$
and column $j$, its label is $j-i$.
The representation theoretical interpretation of these
labels is as follows.
The number associated with the edge linking $\lambda\in\Pi(m)$
and $\mu\in\Pi(m-1)$ is the eigenvalue of the eigenspace
$S^\mu\subset S^\lambda$ of the Jucys-Murphy operator $L_m$,
defined by
$$
L_m=\sum_{i=1}^{m-1} (i,k).
$$

\subsection{Modular representations of $\SG_m$}

On taking the entries of all of the matrices modulo a prime $n$,
the Specht module representations remain well-defined and
we obtain representations of $\SG_m$ over a field of
characteristic $n$.
However, it turns out that these representations are no longer
irreducible in general.
(In fact, they can be both reducible and indecomposable).

\goodbreak
\begin{theorem} {\rm\cite{Ja}}
\begin{enumerate}
\item Over a field of characteristic $n$, the inequivalent irreducible
representations of $\SG_m$ are indexed by
$\Pi_n(m)$, the set of $n$-regular partitions of $m$.
\item For each $\lambda\in\Pi_n(m)$, the corresponding irreducible module
$D^\lambda$ is obtained by $D^\lambda=S^\lambda/{\rm rad}\,S^\lambda$,
where ${\rm rad}\,S^\lambda$ is the maximal proper submodule of $S^\lambda$.
\end{enumerate}
\end{theorem}

\noindent
Despite this result, much remains unknown about the $D^\lambda$.
Of great interest is the search for a combinatorial
algorithm giving the multiplicities of each $D^\mu$ amongst
the composition factors of $S^\lambda$, not least because
this would provide a hitherto unknown direct means to calculate
the dimensions of the $D^\mu$.
The multiplicities of the composition factors of $S^\lambda$
are conveniently written in a matrix with the rows labelled
by $\lambda$ and the columns by $\mu$.
This is the decomposition matrix of $\SG_m$ in characteristic n.
In the case when $m=5$ and $n=2$, we have:
\begin{equation}\label{SGmatrix}
\vcenter{\def\c{\cdot}\halign{\strut\tabskip=5mm $#$\hfil&&\hfil$#$\hfil\cr
&(5) &(41) &(32)\cr
\noalign{\vskip1mm}
(5)  &1&\c&\c\cr
(41) &\c&1&\c\cr
(32) &1&\c&1\cr
(31^2) &2&\c&1\cr
(2^21) &1&\c&1\cr
(21^3) &\c&1&\c\cr
(1^5)  &1&\c&\c\cr}}
\end{equation}

\subsection{The Jantzen-Seitz problem for $\SG_m$}

A modulo $n$ analogue of Theorem \ref{SpechtReduction},
giving a description of the restricted modules
$D^\lambda\!\downarrow^{\SG_m}_{\SG_{m-1}}$, is also currently unknown.

In the $n=2$ case, Benson \cite{Be} considered the somewhat simpler
problem of determining which $\SG_m$-modules $D^\lambda$ remain
irreducible on restriction to $\SG_{m-1}$.
He conjectured that these modules are labelled by the partitions
$\lambda = (\lambda_1,\ldots ,\lambda_r)$ for which the $\lambda_i$
are pairwise distinct and all congruent modulo 2.
Note that in the
characteristic 0 case, Theorem \ref{SpechtReduction}
implies that this problem is solved by the
rectangular partitions $\lambda=(k^l)$,
since it is only in such a case that there is only one way to
remove a node from $\lambda$ to leave a valid partition.

Later, Jantzen and Seitz considered the generalisation
of Benson's conjecture to odd primes $n$. 
They proved that

\begin{theorem} {\rm\cite{JS}}\label{JSTH}
Let $\lambda=(\lambda_1^{a_1},\lambda_2^{a_2},\ldots,\lambda_r^{a_r})$,
where $\lambda_1>\lambda_2>\cdots>\lambda_r>0$ and
$0<a_i<n$ for $i=1,2,\ldots,r$.
Then, over a field of characteristic $n$,
$D^\lambda\!\downarrow^{\SG_m}_{\SG_{m-1}}$ is irreducible if
$$
a_i+\lambda_i-\lambda_{i+1}+a_{i+1}=0 \pmod n,
$$
for $i=1,2,\ldots,r-1$.
\end{theorem}

\noindent
Jantzen and Seitz also conjectured that only for the $\lambda$ specified
here, do the modules $D^\lambda$ remain irreducible on restriction.
This conjecture was proved recently by Kleshchev \cite{Kl1}.

In fact, Kleshchev managed to describe the {\it socle}
({\it i.e.} the sum of all simple submodules) of 
$D^\lambda\!\downarrow^{\SG_m}_{\SG_{m-1}}$.

\begin{theorem} {\rm \cite{Kl3}}
Let $\lambda$ be a $n$-regular partition of $m$. Then
$${\rm Socle}(D^\lambda\!\downarrow^{\SG_m}_{\SG_{m-1}}) \cong \bigoplus_\mu D(\mu)$$
where the sum is over all predecessors $\mu$ of $\lambda$ in Kleshchev's $n$-good
lattice. In particular the sum is multiplicity-free and its description 
is combinatorial.
\end{theorem}

Here Kleshchev's $n$-good lattice is a certain $n$-modular analogue
of Young's lattice. 
It was observed in \cite{LLT1,LLT} that this lattice is the same
as the crystal graph of the basic representation of $\slchap_n$
computed in \cite{MM} (to be described below).

\subsection{The Hecke algebras of type $A$}\label{Hsection}

The (Iwahori-)Hecke algebra $H_m(v)$ of type $A_{m-1}$, 
is the $\C(v)$-algebra generated by the elements 
$T_1,\ldots,T_{m-1}$ subject to the relations:
\begin{eqnarray*}
&&T_iT_{i+1}T_i=T_{i+1}T_iT_{i+1};\\
&&T_iT_j=T_jT_i \qquad \vert{i-j}\vert>1;\\
&&T_i^2=(v-1)T_i+v.
\end{eqnarray*}
If $w\in\SG_m$ and $w=s_{i_1}s_{i_2}\cdots s_{i_l}$ is an
expression for $w$ in terms of a minimal number of generators
of $\SG_m$, define the element $T_w\in H_m(v)$ by
$T_w=T_{i_1}T_{i_2}\cdots T_{i_l}$ (the defining relations
of $H_m(v)$ ensure that $T_w$ is well-defined).

For generic values of $v$, $H_m(v)$ is isomorphic to the group 
algebra of the symmetric group $\SG_m$. 
Thus, it is semisimple, and its 
irreducible representations are parametrised by partitions of $m$.
Again a Specht module construction of the representation matrices
may be undertaken (see Appendix A or \cite{DKLLST}). 

The entries in the representation matrices of the $T_w$ are
now elements of $\Z[v]$.
The branching graph of the Specht modules of $H_m(v)$ is 
the same as for symmetric groups in characteristic 0.
In this case, the number $c$ associated with
an edge linking $\lambda\in\Pi(m)$ and $\mu\in\Pi(m-1)$ indicates
that the eigenvalue of the $v$-Jucys-Murphy
operator $L_m(v)$ is $(v^c-1)/(v-1)$ on the eigenspace
$S^\mu\subset S^\lambda$. 
Here $L_m(v)$ is defined by
$$
L_m(v)=\sum_{i=1}^{m-1} v^{i-m} T_{(i,m)}.
$$

\subsection{Modular representations of Hecke algebras}

In the non-generic case when $v$ is a root of unity, although the Specht
modules remain well defined, they are no longer irreducible in general.

\begin{theorem} {\rm \cite{DJ1}}
Let $v$ be a primitive $n$th root of unity, that is,
$v^n=1$ and $v^k\ne1$ for $1\le k<n$.
\begin{enumerate}
\item The inequivalent irreducible representations of $H_m(v)$ are
indexed by $\Pi_n(m)$, the set of $n$-regular partitions of $m$.
\item For each $\lambda\in\Pi_n(m)$, the corresponding irreducible module
$D^\lambda$ is obtained by $D^\lambda=S^\lambda/{\rm rad}\,S^\lambda$.
\end{enumerate}
\end{theorem}

\noindent
Although the same labels $D^\lambda$ are traditionally used for
the irreducible representations of $\SG_m$ over a field of
characteristic $n$, and for the irreducible representations of
$H_m(\sqrt[n]{1})$,
the representations themselves should not be
identified in any sense. Indeed, they may even have different
dimensions. Moreover those of $\SG_m$ are defined only for
$n$ prime.

As in the case of modular representations of $\SG_m$,
one defines 
the decomposition matrices of $H_m(\sqrt[n]{1})$.
Apart from their independent interest, they may be viewed as a
stepping stone to the modular case of $\SG_m$.
Compare the following table for $H_5(-1)$ (here $m=5$ and $n=2$)
with that given in Section \ref{SGsection}:
\begin{equation}\label{Hmatrix}
\vcenter{\def\c{\cdot}\halign{\strut\tabskip=5mm $#$\hfil&&\hfil$#$\hfil\cr
&(5) &(41) &(32)\cr
\noalign{\vskip1mm}
(5)  &1&\c&\c\cr
(41) &\c&1&\c\cr
(32) &\c&\c&1\cr
(31^2) &1&\c&1\cr
(2^21) &\c&\c&1\cr
(21^3) &\c&1&\c\cr
(1^5)  &1&\c&\c\cr}}
\end{equation}

\noindent
Very recently, a combinatorial  algorithm for calculating the
decomposition matrices for $H_m(\sqrt[n]{1})$ was conjectured 
in \cite{LLT1,LLT}. This algorithm computes in fact Kashiwara's 
global crystal basis of the basic representation of the quantum 
affine algebra $U_q(\slchap_n)$ (see \ref{LGCBsection}).

Lascoux, Leclerc and Thibon's interpretation of these bases in terms of 
the Grothendieck ring of representations of $H_m(\sqrt[n]{1})$ has now 
been proved by Ariki and Grojnowski.

This machinery enables us to determine for which $\lambda$ the restricted 
module $D^\lambda\!\downarrow^{H_m(\sqrt[n]{1})}_{H_{m-1}(\sqrt[n]{1})}$ 
is irreducible, and shows how to relate this Jantzen-Seitz type problem 
for Hecke algebras with the decomposition of tensor products of level 1 
representations of $\slchap_n$.



\def\smbox#1{\smyoungd{
  \multispan3\hrulefill\cr &#1&\cr\multispan3\hrulefill\cr}}

\section{The Fock representation of $U_q(\slchap_n)$}

For our purposes, the relevant realisation of $\slchap_n$
and $U_q(\slchap_n)$ is the Fock representation.
Moreover, the reduction from $H_m(v)$ to $H_m(\sqrt[n]{1})$
is reflected in the reduction from $\glchap_\infty$
to $\slchap_n$, so we shall first recall the Fock 
representation of $\glchap_\infty$.

\subsection{The infinite rank affine algebra $\glchap_\infty$}

Let $\overline{\goth gl}_\infty$ denote the Lie
algebra of complex $\Z\times \Z$-matrices
$A=(a_{ij})$ such that $a_{ij}=0$ for $|i-j|\gg0$. 
Therefore, $A$ has only a finite number of nonzero diagonals,
and matrix multiplication makes sense, which allows
to define the Lie bracket as the commutator $AB-BA$.
Introduce the following block decomposition  
$$A_{++}=(a_{ij})_{i,j>0},\ 
A_{+-}=(a_{ij})_{i>0,j\le 0},\
A_{-+}=(a_{ij})_{i\le 0,j>0},\
A_{--}=(a_{ij})_{i,j\le 0}.
$$
Then (see \cite{JM}) $\glchap_\infty$ is defined as the one-dimensional central 
extension 
$$\glchap_\infty = \overline{\goth gl}_\infty \oplus \C\,c$$
with Lie bracket
\begin{eqnarray*}
[A,B]&=&AB-BA + \tr(A_{-+}B_{+-}-A_{+-}B_{-+})\,c, \quad 
(A,B\in \overline{\goth gl}_\infty), \\[2mm]
[c,\glchap_\infty]& = & 0. 
\end{eqnarray*}
The Chevalley generators of $\glchap_\infty$
are expressed in terms of the matrix units $E_{ij}$ by
\begin{eqnarray*}
e^\infty_i&=&E_{i,i+1},\\[2mm]
f^\infty_i&=&E_{i+1,i}, \\[2mm]
h^\infty_i&=&[e_i^\infty,f_i^\infty]=E_{i,i}-E_{i+1,i+1}+\delta_{i0}c,
\quad
(i\in\Z).
\end{eqnarray*}

By \cite{JM}, $\glchap_\infty$ acts in a natural way on 
the polynomial ring 
$$\F=\C[x_i, i\ge 1],$$ 
where the variable $x_i$ has degree $i$.
This graded space is called the (bosonic) Fock space. 
The action is best described on a distinguished basis of 
homogeneous elements indexed by partitions,
the Schur polynomials
$$
s_\lambda = \sum_{\mu = (1^{m_1}\ldots r^{m_r})}
\chi_\lambda(\mu)\,{x_1^{m_1}\cdots x_r^{m_r}\over m_1!\cdots m_r!} ,
$$
where for $\lambda$, $\mu$ partitions of $m$,
$\chi_\lambda(\mu)$ denotes the irreducible character
$\chi_\lambda$ of $\SG_m$ evaluated on a element of cycle-type $\mu$.
 
To write down the action, we need some combinatorial definitions.
Let $\gamma$ be the node of the partition $\lambda$
situated on row $i$ and column $j$.
The content $c=c(\gamma)$ of $\gamma$ is defined by $c(\gamma)=j-i$.
The node $\gamma$ is then said to be a $c$-node.
A removable $c$-node of $\lambda$ is a $c$-node $\gamma$ within
$\lambda$, which if removed from $\lambda$ leaves a valid partition $\mu$.
In this case, we also say that $\gamma$ is an addable $c$-node
of $\mu$ and we use freely the notation
$\lambda/\mu=\smbox{c}$
or $\lambda/\,\smbox{c}=\mu$
or $\mu\cup\smbox{c}=\lambda$.

\begin{example}{\rm
To illustrate these definitions,
consider the partition $(7,5,5,2,1,1)$, which has Young diagram:
$$
\smyoungd{
\multispan{15}\hrulefill\cr &0&&1&&2&&3&&4&&5&&6&\cr
\multispan{15}\hrulefill\cr &\bar1&&0&&1&&2&&3&\cr
\multispan{11}\hrulefill\cr &\bar2&&\bar1&&0&&1&&2&\cr
\multispan{11}\hrulefill\cr &\bar3&&\bar2&\cr
\multispan{5}\hrulefill\cr  &\bar4&\cr
\multispan{3}\hrulefill\cr  &\bar5&\cr
\multispan{3}\hrulefill\cr}.
$$
Here each node has been filled with its content,
and for convenience
the negative entries have the minus sign above the digit.
We see that $\lambda$ has a removable $c$-node only
if $c=-5,-2,2,6$ and an addable $c$-node only if $c=-6,-3,-1,4,7$.
Taking an example from each of these cases, we have
$\lambda/\mu=\smbox{\bar2}$ for $\mu=(7,5,5,1,1,1)$ and
$\nu/\lambda=\smbox{7}$ for $\nu=(8,5,5,2,1,1)$.}
\end{example}

We can now state
\begin{theorem} {\rm \cite{JM}}\label{glchapRep}
The Lie algebra $\glchap_\infty$ acts on $\F$ as follows:
\begin{eqnarray}
e^\infty_i s_\lambda
&=&\cases{ s_\mu\phantom{-{}} &if $\lambda/\mu=\smbox{i}$\,;\cr
           0 &otherwise,\cr} \\
f^\infty_i s_\lambda
&=&\cases{ s_\nu\phantom{-{}} &if $\nu/\lambda=\smbox{i}$\,;\cr
           0 &otherwise,\cr} \\
h^\infty_i s_\lambda
&=&\cases{ -s_\lambda &if $\lambda$ has a removable $i$-node;\cr
           s_\lambda &if $\lambda$ has an addable $i$-node;\cr
           0 &otherwise,\cr} \\
c \,s_\lambda & = & s_\lambda .
\end{eqnarray}
This is an irreducible highest weight representation with
highest weight $\Lambda_0$, that is, the highest weight
vector $1=s_0$ satisfies
$$
h_i^\infty s_0 = \delta_{i0}\, s_0 .
$$
\end{theorem}

\noindent
Up to partitions with 5 nodes, the action of the generators
is given by the following graph.

\medskip
\centerline{\epsfbox{young_fock.eps}}
\medskip
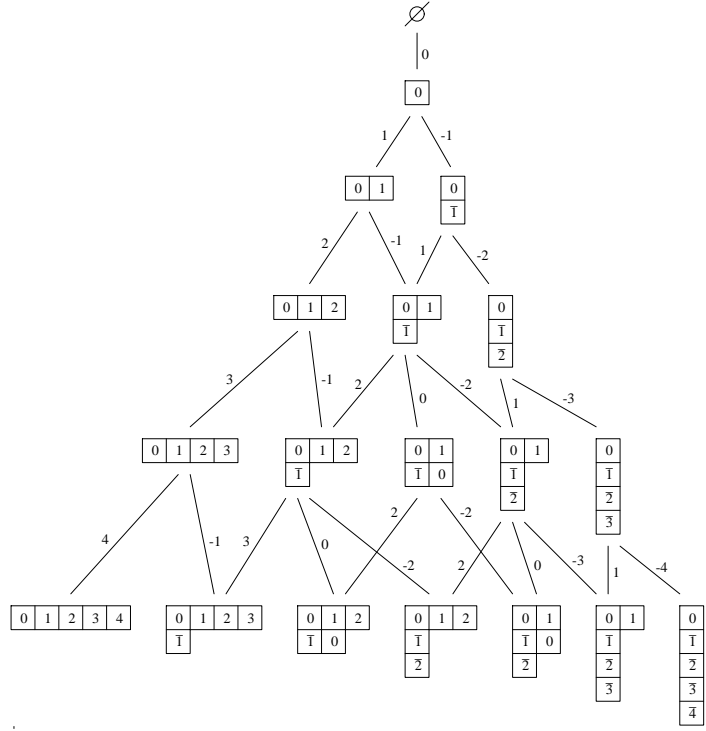
\centerline{Figure 2.}
\medskip

\noindent We see that this graph is identical (apart from the
inclusion of the empty partition $\emptyset$) to that in Fig. 1.

Moreover, if we define the complete restriction operator
$$
e^\infty=\sum_{i=-\infty}^\infty e^\infty_i,
$$
the action of $e^\infty$ on $s_\lambda$ gives
$$
e^\infty s_\lambda = \sum_{\mu\in{\cal R}(\lambda)} s_\mu,
$$
where ${\cal R}(\lambda)$ is as given in Theorem \ref{SpechtReduction}.
Thus, the Fock space representation of $\glchap_\infty$ yields
the branching rule for $\SG_m$ in characteristic 0
(or $H_m(v)$ for generic~$v$).

In what follows, a similar correspondence 
connecting the Hecke algebras $H_m(\sqrt[n]{1})$ and
the Fock space representation of $\slchap_n$
is explained.

\subsection{The embedding of $\slchap_n$ into $\glchap_\infty$}\label{Sec4.2}

It is known that many affine Lie algebras
can be realized as subalgebras of $\glchap_\infty$.
This type of embedding has been used by the Kyoto group
to exhibit relations between different hierarchies 
of soliton equations (see \cite{JM}). 
In particular, $\slchap_n$ can be realized
as a subalgebra of $\glchap_\infty$ by taking as Chevalley generators
of the former
\begin{eqnarray*}
e_i&=& \sum_{k=i\,(\mod n)} e^\infty_k,\\
f_i&=& \sum_{k=i\,(\mod n)} f^\infty_k,\\
h_i&=& \sum_{k=i\,(\mod n)} h^\infty_k,\\
D&=& -\sum_{k\in \Z} \left[{k\over n}\right]\, E_{ii}.
\end{eqnarray*}
Restricting to $\slchap_n$ the action of $\glchap_\infty$
on $\F$, one gets the Fock space representation of $\slchap_n$.
Clearly, one has
\begin{eqnarray*}
\displaystyle f_i s_\lambda&=&
  \sum_{\mu:\mu/\lambda=\smbox{i}} \, s_\mu,\\
\displaystyle e_i s_\lambda&=&
  \sum_{\mu:\lambda/\mu=\smbox{i}} \, s_\mu.
\end{eqnarray*}
This action may be easily obtained from the graph in Fig.~2,
after each of the labels
(both in the partitions and on the edges) is replaced by 
its corresponding modulo $n$ value \cite{DJKMO}.

Under the action of $\slchap_n$, the Fock space $\F$ 
decomposes into
\begin{equation}\label{DECOMPF}
\F = {\cal T}_n \oplus {\cal S}_n
\end{equation}
where 
$${\cal T}_n=U(\slchap_n)\,s_0 = \C[x_i, i\not\equiv 0 {\ \rm mod\ } n],$$
and ${\cal S}_n$ is the ideal generated by the $x_{kn}$,
$(k>0)$.
Moreover, as $\slchap_n$-modules,
$$
{\cal T}_n \cong \F/{\cal S}_n \cong V(\Lambda_0),
$$
the basic representation.
Note that the natural map from $\F$ to $\F/{\cal S}_n$ obtained by
setting $x_{kn} = 0$ is precisely the reduction of the KP-hierarchy
to the ${\rm KP}_n$-hierarchy (where ${\rm KP}_2 = {\rm KdV}$)
(see \cite{JM}).

In our setting, the same map becomes the decomposition
map from the Grothendieck group of $H_m(v)$-modules
to that of $H_m(\sqrt[n]{1})$-modules (for all $m$).

\subsection{The quantum affine algebra $U_q(\slchap_n)$}

For $k,m\in\Z$, we define the following $q$-integers, factorials,
and binomials

$$
[k]_q = {q^k- q^{-k} \over q - q^{-1}}, 
\quad
[k]_q! = [k]_q\,[k-1]_q\,\cdots [1]_q, 
\quad
\hbox{and}
\quad
{m \brack k}_q = {[m]_q!\over [m-k]_q!\,[k]_q!} \,.
$$

The quantum affine algebra $U_q(\slchap_n)$ is
the unital associative algebra over $\Q(q)$
generated by the symbols $e_i,\ f_i,\ 0\le i \le n-1,$ and
$q^h$ for $h\in P^\vee$, subject to the relations:
$$
q^h\,q^{h'}  =   q^{h+h'};  \quad  q^0 = 1;
$$ $$
q^h e_j q^{-h} = q^{\<\alpha_j,h\>} e_j ;
$$ $$
q^h f_j q^{-h} = q^{-\<\alpha_j,h\>} f_j ;
$$ $$
{}[e_i,f_j] = \delta_{ij} {q^{h_i}- q^{-h_i} \over q - q^{-1}} ;
$$ $$
\sum_{k=0}^{1-\<\alpha_i,h_j\>} (-1)^k 
 {1-\<\alpha_i,h_j\> \brack k}_q 
 e_i^{1-\<\alpha_i,h_j\> -k} e_j e_i^k = 0  \quad (i\not = j) ;
$$ $$
\sum_{k=0}^{1-\<\alpha_i,h_j\>} (-1)^k
 {1-\<\alpha_i,h_j\> \brack k}_q
 f_i^{1-\<\alpha_i,h_j\> -k} f_j f_i^k = 0  \quad (i\not = j) .
$$

The definitions concerning $U_q(\slchap_n)$-modules are direct analogues 
of those given for $\slchap_n$. Let $M$ be a $U_q(\slchap_n)$-module.
For each weight $\Lambda \in P$, the subspace of $M$ defined by
$$
M_\Lambda = 
\{ v\in M \ | \ q^h\,v = q^{\<\Lambda,h\>} \,v, \ h\in P^\vee \}
$$
is called the weight space of weight $\Lambda$ of $M$.
If $v\in M_\Lambda$ then we write $\wt(v)=\Lambda$, and
call such $v$ a weight vector of weight $\Lambda$.

The module $M$ is said to be integrable if:
\begin{enumerate}
\item $M=\bigoplus_{\Lambda \in P} M_\Lambda$;
\item ${\rm dim}\,M_\Lambda < \infty$ for each $\Lambda \in P$;
\item for each $i= 0,1,\ldots ,n-1,$ $M$ decomposes into a direct sum
of finite dimensional $U_i$-modules, where $U_i$ denotes the subalgebra 
of $U_q(\slchap_n)$ generated by $e_i$, $f_i$, $q^{h_i}$ and $q^{-h_i}$.
\end{enumerate}

A weight vector $v\in M$ for which $e_iv=0$ for all
$i=0,1,\ldots,n-1$, is said to be a highest weight vector.
If there exists a highest weight vector $v\in M$ such that
$M = U_q(\slchap_n)\, v$, then $M$ is said to be a highest weight module.
The weight of $v$ is called the highest weight of $M$.

For $l\ge0$ and all $\Lambda\in P^+_l$, there exists
(up to equivalence)
a unique integrable highest weight module $V(\Lambda)$ 
with highest weight $\Lambda$,
and it is irreducible \cite{ChariPressleyBook}.

\subsection{The Fock space representation of $U_q(\slchap_n)$}

We require further notation concerning removable and addable nodes.
From now on, we say that $\gamma$ is a $i$-node
if its content $c$ is such that $c=i\pmod n$. 
\begin{definition}
Let $\lambda$ be a partition, $0\le i<n$, $\gamma$ be an
addable $i$-node of $\lambda$, and $\nu=\lambda\cup\gamma$.
We set
\begin{enumerate} 
\item
$N_i(\lambda) = \#\{$ addable $i$-nodes of $\lambda\ \}
   - \#\{$ removable $i$-nodes of $\lambda\ \}$,
\item
$N_i^l(\lambda,\nu) = \#\{$ addable $i$-nodes of $\lambda$ situated 
   strictly to the {\it left} of $\gamma\ \}$
   \newline\hbox to25mm{\hfil}
   $-\#\{$ removable $i$-nodes of $\lambda$ situated 
   to the {\it left} of $\gamma\ \}$,
\item
$N_i^r(\lambda,\nu) = \#\{$ addable $i$-nodes of $\lambda$ situated
   strictly to the {\it right} of $\gamma\ \}$
   \newline\hbox to25mm{\hfil}
   $-\#\{$ removable $i$-nodes of $\lambda$ situated
   to the {\it right} of $\gamma\ \}$,
\item
$N^0(\lambda) = \#\{$ 0-nodes of $\lambda\ \}$.
\end{enumerate}
\end{definition}

\begin{center}
\unitlength=0.0125in
\begin{picture}(255,125)(0,-10)
\path(0,0)(20,0)(20,10)
	(40,10)(40,30)(90,30)
	(90,60)(110,60)(110,70)
	(120,70)(120,80)(150,80)
	(150,100)(180,100)(180,110)
	(0,110)(0,0)
\path(80,30)(80,40)(90,40)
\path(57.071,5.757)(50.000,10.000)(54.243,2.929)
\path(50,10)(60,0)(170,0)
\path(137.071,65.757)(130.000,70.000)(134.243,62.929)
\path(130,70)(140,60)(255,60)
\put(65,5){\makebox(0,0)[lb]{\raisebox{0pt}[0pt][0pt]{\shortstack[l]{{\tenrm nodes to the left of $\gamma$}}}}}
\put(145,65){\makebox(0,0)[lb]{\raisebox{0pt}[0pt][0pt]{\shortstack[l]{{\tenrm nodes to the right of $\gamma$}}}}}
\put(82,34){\makebox(0,0)[lb]{\raisebox{0pt}[0pt][0pt]{\shortstack[l]{{\tenrm   $\gamma$}}}}}
\end{picture}
\end{center}
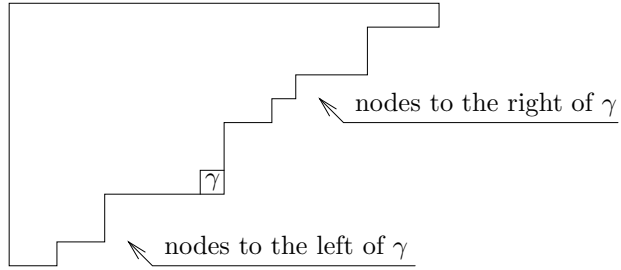
\centerline{Figure 3.}
\medskip

The following construction of the Fock space representation
of $U_q(\slchap_n)$
is due to Hayashi \cite{Hay}.
We give a formulation by Misra and Miwa \cite{MM}.

\begin{theorem}[$q$-Fock space]
{\rm \cite{Hay,MM}} Let the $q$-Fock space be defined by
$$
\F = \bigoplus_{\lambda\in\Pi} \Q(q)\,v_\lambda.
$$
Then $U_q(\slchap_n)$ acts on $\F$ by
\begin{eqnarray*}
q^{h_i}\,v_\lambda&=&
  q^{N_i(\lambda)}\,v_\lambda\,,\\
q^D\,v_\lambda&=&
  q^{-N^0(\lambda)}\,v_\lambda\,,\\
\displaystyle f_i v_\lambda&=&
  \sum_{\nu:\nu/\lambda=\smbox{i}} q^{N_i^r(\lambda,\nu)} \, v_\nu\,,\\
\displaystyle e_i v_\lambda&=&
  \sum_{\nu:\lambda/\nu=\smbox{i}} q^{-N_i^l(\nu,\lambda)} \, v_\nu\,.
\end{eqnarray*}
This representation is integrable.
\end{theorem}

The \lq\lq vacuum\rq\rq\ vector $v_\emptyset$ labelled by the
empty partition, is a highest weight vector of weight $\Lambda_0$.
The decomposition of $\F$ under the action 
of $U_q(\slchap_n)$ is similar to (\ref{DECOMPF}).
In particular one has
$$
U_q(\slchap_n)\,v_\emptyset \cong V(\Lambda_0).
$$

In the following section, we will describe a basis for
this subrepresentation of $\F$. Of interest to us will 
be the transition matrices between
that basis and the basis $(v_\lambda)$ of $\F$.

\section{Canonical bases}

Canonical bases of the irreducible highest weight modules of
affine algebras have been introduced by Kashiwara as classical
limits ($q\mapsto 1$) of the so-called global crystal bases of 
the corresponding quantum affine algebras. 
The first step towards obtaining the global crystal basis is to
determine the crystal basis at $q=0$ and the associated crystal graph, 
a combinatorial object describing the action of the quantum algebra
in the crystal limit $q\mapsto 0$. 

\subsection{The crystal basis at $q=0$}

In this section, we will describe following \cite{MM} 
the crystal graph
of $\F$ and isolate in it the crystal graph of
$V(\Lambda_0)$.

To define a crystal basis at $q=0$ (or simply a crystal basis),
one needs to introduce Kashiwara's operators.

On restricting to $U_i$, the module $M$
decomposes into a direct sum of simple $U_i$-modules, with
each $(l+1)$-dimensional $U_i$-module isomorphic to the
$U_q(\Sl_2)$-module $V_l$ having highest weight $l\Lambda_1$:
$$
M \cong \bigoplus_k V_{l_k}
$$
Now, for arbitrary $l$, write
$$
V_l = \bigoplus_{k=0}^l \Q(q)\,u_k^{(l)}\,,
$$
where
$$
e_i\,u_0^{(l)} = f_i u_l^{(l)} = 0
\quad{\rm and}\quad
u_k^{(l)} = {f_i^k\over [k]_q!}\,u_0^{(l)}\,.
$$
Then the action of $\tilde{e_i}$ and $\tilde{f_i}$ on $V_l$
is defined by:
$$
\tilde{f_i} \,u_k^{(l)} = u_{k+1}^{(l)}\,, \quad
\tilde{e_i} \,u_k^{(l)} = u_{k-1}^{(l)}\,,
$$
where we understand $u_{-1}^{(l)} = u_{l+1}^{(l)} = 0$.
The endomorphisms $\tilde{e_i},\ \tilde{f_i}$ are then extended to
$M$ using the above isomorphism of $U_i$-modules
between $M$ and $\bigoplus_k V_{l_k}$.

Thus, any integrable $U_q(\slchap_n)$-module $M$ is equipped with
a family of endomorphisms $\tilde{e_i}$ and $\tilde{f_i}$.

Let $A\subset \Q(q)$ denote the ring of rational functions without
a pole at $q=0$.
A {\it crystal lattice} of $M$ is a free $A$-module $L$ such that
$M=\Q(q)\otimes_A L,\ L = \oplus_\Lambda L_\Lambda$
where $L_\Lambda = L \cap M_\Lambda$,
and 
$$ 
\tilde{e_i} L \subset L, \quad
\tilde{f_i} L \subset L, \quad
(i=0,1,\ldots n-1).
$$
In other words, $L$ spans $M$ over $\Q(q)$, $L$ is compatible
with the weight space decomposition of $M$ and is stable
under Kashiwara's operators.
It follows that $\tilde e_i, \, \tilde f_i$ induce endomorphisms
of the $\Q$-vector space $L/qL$ that we shall still denote by 
$\tilde e_i, \, \tilde f_i$.
Now Kashiwara defines a {\it crystal basis}
 of $M$ (at $q=0$) to be a pair
$(L,B)$ where $L$ is a crystal lattice in $M$ and $B$ is a basis of $L/qL$ such
that $B=\sqcup B_\Lambda$ where $B_\Lambda = B \cap (L_\Lambda/qL_\Lambda)$, and
$$
\tilde e_i B \subset B\sqcup \{0\}\,, \quad
\tilde f_i B \subset B\sqcup \{0\}\,, \quad i = 0,\ldots ,n-1\,,
$$
$$
\tilde e_i v = u \quad \Longleftrightarrow \quad \tilde f_i u = v ,\, \quad u,v \in B,
\,
\quad i = 0,\ldots ,n-1\,.
$$
Kashiwara has proven the following existence and uniqueness result for crystal
bases \cite{Ka1,Ka2}.
\begin{theorem}\label{EXIST}
Any integrable $U_q(\slchap_n)$-module $M$ has a crystal basis $(L,B)$. Moreover,
if $(L',B')$ is another crystal basis of $M$, then there exists a
$U_q(\slchap_n)$-automorphism of $M$ sending $L$ to $L'$ thence inducing an isomorphism
of vector spaces from $L/qL$ to $L'/qL'$ which sends $B$ to $B'$. In particular,
if $M= V(\Lambda)$ is irreducible, its crystal basis $(L(\Lambda),B(\Lambda))$ is
unique up to an overall scalar multiple. It is given by
\begin{equation}
L(\Lambda)= \sum_{0\le i_1,i_2,\ldots,i_r\le n-1}  A
\, \tilde f_{i_1}\tilde f_{i_2}\cdots \tilde f_{i_r}\, u_\Lambda ,
\end{equation}
\begin{equation}
B(\Lambda) = \{
\tilde f_{i_1}\tilde f_{i_2}\cdots \tilde f_{i_r}\, u_\Lambda\ \mod qL(\Lambda)\, |\,
0\le i_1,\ldots,i_r\le n-1\}\backslash\{0\} ,
\end{equation}
where $u_\Lambda$ is a highest weight vector of $V(\Lambda)$.
\end{theorem}

It follows that to each integrable $U_q(\slchap_n)$-module $M$, one can associate
a well-defined coloured graph $\Gamma(M)$ whose vertices are labelled by the
elements of $B$ and whose edges describe the action of the operators $\tilde f_i$ :
$$
u \stackrel{i}{\longrightarrow} v \ \Longleftrightarrow \ \tilde f_iu = v\,.
$$
$\Gamma(M)$ is called the {\it crystal graph} of $M$.

A crystal basis and the crystal graph of the $q$-Fock space $\F$
have been determined by Misra and Miwa.

\begin{theorem}\label{FockLcbTheorem} {\rm \cite{MM}}
The pair $(L,B)$ where
$$
L=\bigoplus_{\lambda\in\Pi} A\,v_\lambda\,,
$$
and $B$ is the basis of $L/qL$ given by
$$
B=\{v_\lambda\mod qL\mid\lambda\in\Pi\}\,.
$$
is a crystal basis of $\F$.
\end{theorem}

To describe the crystal graph we introduce the notion of {\it good node}.
Let $\lambda$ be a partition and $0\le i<n$.
Construct a sequence of $A$s and $R$s, by
from left to right, scanning the columns of $\lambda$ and noting
the presence of an addable $i$-node by $A$ and the presence of
a removable $i$-node by $R$.
It will be useful to attach a subscript to each of these symbols
to indicate the column from which it arose.

Now recursively remove $RA$ pairs (together with their subscripts)
from this sequence until none remain.
The sequence will then be of the form 
$$AA\cdots AR\cdots RRR.$$
The node corresponding to the leftmost $R$ is termed a
good (removable) $i$-node and that corresponding to the
rightmost $A$ is termed a good addable $i$-node (note that
there is at most one of each).

\begin{example} {\rm To illustrate this definition, consider $n=3$ and
the partition $\lambda=(16,13,11,10,9,8,7,5,2)$:
$$
\smyoungd{
\multispan{33}\hrulefill\cr
&0&&1&&2&&0&&1&&2&&0&&1&&2&&0&&1&&2&&0&&1&&2&&0&\cr
\multispan{33}\hrulefill\cr
&2&&0&&1&&2&&0&&1&&2&&0&&1&&2&&0&&1&&2&\cr
\multispan{27}\hrulefill\cr
&1&&2&&0&&1&&2&&0&&1&&2&&0&&1&&2&\cr
\multispan{23}\hrulefill\cr
&0&&1&&2&&0&&1&&2&&0&&1&&2&&0&\cr
\multispan{21}\hrulefill\cr
&2&&0&&1&&2&&0&&1&&2&&0&&1&\cr
\multispan{19}\hrulefill\cr
&1&&2&&0&&1&&2&&0&&1&&2&\cr
\multispan{17}\hrulefill\cr
&0&&1&&2&&0&&1&&2&&0&\cr
\multispan{15}\hrulefill\cr
&2&&0&&1&&2&&0&\cr
\multispan{11}\hrulefill\cr
&1&&2&\cr
\multispan{5}\hrulefill\cr}\;.
$$
Here there are addable $0$-nodes in columns 1,3,9,12 and 14,
and removable $0$-nodes in columns 5,7,10 and 16.
Thus we form the following sequence:
$$
A_1 A_3 R_5 R_7 A_9 R_{10} A_{12} A_{14} R_{16}.
$$
The removal procedure first disposes of $R_7 A_9$ and $R_{10} A_{12}$,
and then disposes of $R_5 A_{14}$ so that:
$$
A_1 A_3 R_{16}
$$
remains. Therefore $\lambda$ has a good addable node in the 3rd column
and a good removable node in the 16th column.
In the case of 1-nodes and 2-nodes, we first obtain the sequences:
\begin{eqnarray*}
&&A_6 A_8 R_9 A_{11} A_{17}\\
\noalign{\noindent\hbox{and}}
&&R_2 R_8 A_{10} R_{11} R_{13},
\end{eqnarray*}
respectively which, after the removal procedure, produce:
\begin{eqnarray*}
&&A_6 A_8 A_{17}\\
\noalign{\noindent\hbox{and}}
&&R_2 R_{11} R_{13}
\end{eqnarray*}
respectively.
Thus $\lambda$ has a good addable 1-node in column 17, but no good
removable 1-node, and a good removable 2-node in column 2, but
no good addable 2-node.}
\end{example}

\begin{theorem}{\rm\cite{MM}}
The crystal graph $\Gamma_n$ of $\F$ is the graph with vertices
labelled by $\Pi$, the set of partitions, and
whose arrows are given by:
$$ 
\lambda \stackrel{i}{\longrightarrow} \mu
\Longleftrightarrow
\hbox{$\mu$ is obtained from $\lambda$ by adding a good addable $i$-node.}
$$
\end{theorem}
In such a case, we write $\tilde{f_i}(\lambda)=\mu$ and
$\tilde{e_i}(\mu)=\lambda$.
For each $\lambda\in\Pi$, the largest integer $k$ such that
$\tilde{e_i}^k(\lambda)\ne0$ is denoted $\varepsilon_i(\lambda)$.
The largest integer $k$ such that $\tilde{f_i}^k(\lambda)\ne0$ is
denoted $\varphi_i(\lambda)$.

\begin{example}
{\rm Up to partitions of weight 5,
the crystal graph $\Gamma_2$
is as follows:

\medskip
\centerline{\epsfbox{sl2_0f.eps}}
\medskip
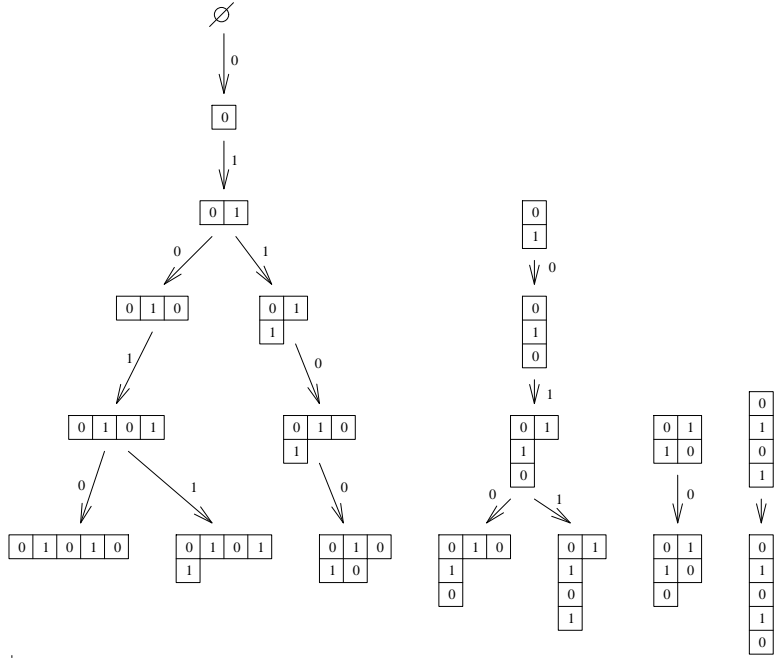
\centerline{Figure 4.}}
\end{example}

It is clear from the definition that the crystal graph $\Gamma(M)$ of the direct
sum $M=M_1\oplus M_2$ of two $U_q(\slchap_n)$-modules is the disjoint union of
$\Gamma(M_1)$ and $\Gamma(M_2)$. It follows from the complete reducibility
of $M$ that the connected components of $\Gamma(M)$ are the crystal graphs 
of the irreducible components of $M$ \cite{Ka1,Ka2}.

Thus, in Fig.~4, the splitting of $\Gamma_2$ into its
connected components reflects the decomposition of the 
Fock representation of $U_q(\slchap_2)$.

In particular, the connected component of $\emptyset$
in $\Gamma_n$ is the crystal graph of the $U_q(\slchap_n)$-module
$V(\Lambda_0)$.
Its set of vertices is $\Pi_n$, the set of $n$-regular partitions.
The crystal graph of $V(\Lambda_0)$ for $n=3$ is shown in
Fig. 5 (up to partitions of weight 8).

It is easily seen that each connected component of $\Gamma_n$
is headed by a unique partition $\mu$
of the form $\lambda^{*n}$, where if
$\lambda=(\lambda_1^{a_1},\lambda_2^{a_2},\ldots,\lambda_r^{a_r})$,
we define
$\lambda^{*n}=(\lambda_1^{na_1},\lambda_2^{na_2},\ldots,\lambda_r^{na_r})$.
This is the combinatorial image of the known decomposition
under $U_q(\slchap_n)$
$$
\F \cong \bigoplus_{k\ge 0} V(\Lambda_0-k\delta)^{\oplus p(k)} \,,
$$
where $p(k)$ denotes the number of partitions of $k$.
\subsection{Tensor products}

One of the main properties of crystal bases at $q=0$ is that they behave 
well under tensor products \cite{Ka1,Ka2}. 
\begin{theorem}\label{TP}
Let $(L_1,B_1)$ and $(L_2,B_2)$ be crystal bases of integrable 
$U_q(\slchap_n)$-modules $M_1$ and $M_2$. 
Let $B_1\otimes B_2$ denote the basis $\{ u\otimes v,\ u\in B_1,\, 
v\in B_2\}$
of $(L_1/qL_1 )\otimes (L_2/qL_2)$ (which is isomorphic to
$(L_1\otimes L_2)/q(L_1\otimes L_2)$). Then, $(L_1\otimes L_2,\,B_1\otimes B_2)$ is a
crystal basis of $M_1\otimes M_2$, the action of $\tilde e_i$, $\tilde f_i$ on
$B_1\otimes B_2$ being given by
\begin{eqnarray*}
\tilde{f_i}(u\otimes v)&=&\cases{
\tilde{f_i}u\otimes v&if $\varphi_i(u)>\varepsilon_i(v)$;\cr
u\otimes\tilde{f_i}v&otherwise,\cr}\\
\tilde{e_i}(u\otimes v)&=&\cases{
\tilde{e_i}u\otimes v&if $\varphi_i(u)\ge\varepsilon_i(v)$;\cr
u\otimes\tilde{e_i}v&otherwise.\cr}
\end{eqnarray*}
\end{theorem}

Taking into account the fact that the decomposition of an integrable 
module $M$ into irreducible components is reflected by the splitting 
of its crystal graph into connected components, one sees that 
Theorem~\ref{TP} gives a powerful way of computing tensor product 
multiplicities.
  
Let $\Lambda^\prime$, $\Lambda^{\prime\prime}$, $\Lambda$ be dominant integral
weights. The multiplicity
$c_{\Lambda^\prime\,\Lambda^{\prime\prime}}^\Lambda$
of $V(\Lambda)$ in the
tensor product $V(\Lambda^\prime)\otimes V(\Lambda^{\prime\prime})$
is equal to the number of vertices $b_1\otimes b_2$ of 
$B(\Lambda^\prime)\otimes B(\Lambda^{\prime\prime})$ that satisfy
$$
\wt(b_1\otimes b_2) = \Lambda
\quad \mbox{and} \quad \tilde{e}_i (b_1\otimes b_2) = 0 \quad (i=0,\ldots ,n-1).
$$
By Theorem~\ref{TP}, this last condition is equivalent to the
fact that $b_1=b_{\Lambda^\prime}$,  the origin of the crystal graph
of $V(\Lambda^\prime)$, and 
$$
\varepsilon_i(b_2)\le \<\Lambda^\prime,h_i\>\qquad i=0,1,\ldots,n-1.
$$
Hence we get,
\begin{corollary}\label{TPcorollary}
The multiplicity $c_{\Lambda^\prime\,\Lambda^{\prime\prime}}^\Lambda$
of $V(\Lambda)$ in 
$V(\Lambda^\prime)\otimes V(\Lambda^{\prime\prime})$
is equal to the number of vertices $b_2$ of $B(\Lambda^{\prime\prime})$
such that
$$
\wt(b_2)=\Lambda - \Lambda^\prime \quad \mbox{  and  } \quad
\varepsilon_i(b_2)\le \<\Lambda^\prime,h_i\>\qquad i=0,1,\ldots,n-1.
$$
\end{corollary}

\subsection{The lower global crystal basis}\label{LGCBsection}

The crystal basis of a $U_q(\slchap_n)$-module $M$
is not a basis of $M$, but only a combinatorial
object reflecting the structure of $M$ in the crystal
limit $q\mapsto 0$.
To get a true basis of $M$, one needs to ``globalize"
this notion (Kashiwara's terminology) and define
a global crystal basis.
Kashiwara's theory enables such a basis to be obtained
for all irreducible highest weight modules $V(\Lambda)$.
However, if $M$ is not irreducible, there is no general
way to obtain a global basis of $M$.
For example, this is the situation for the Fock space $\F$
(note however, \cite{LT}).

In fact there are two global bases of $V(\Lambda)$, the lower one and
the upper one, which are dual one to the other.
In this section we consider the first and
recall the definition of the lower global crystal basis
of $V(\Lambda_0)$.
We need to introduce an involution $v \longrightarrow \overline{v}$
of $V(\Lambda_0)$. We start from the involution $P\longrightarrow \overline{P}$
of $U_q(\slchap_n)$ defined as the ring automorphism satisfying
$$
\overline{q} = q^{-1},\quad  \overline{q^h} = 
q^{-h}, \qquad (h \in P^\vee) \,,
$$
$$\overline{e_i} = e_i,\quad \overline{f_i} = f_i, \quad i=0,1,\ldots ,n-1 \,.$$
Then, for $v=P\,v_\emptyset \in V(\Lambda_0)$, we set
$\overline{v}=\overline{P}\,v_\emptyset$.
Finally, we denote by $U_\Q^-$ the sub-$\Q[q,q^{-1}]$-algebra of
$U_q(\slchap_n)$ generated by $f_i^{(k)}:=f_i^k/[k]!$, and we
set $V_\Q(\Lambda_0)=U_\Q^-\,v_{\emptyset}$.
We can now state
\begin{theorem}\label{LGCBTheorem} {\rm\cite{Ka2}}
There exists a unique $\Q[q,q^{-1}]$-basis
$\{G(\mu)\mid \mu\in \Pi_n\}$ of $V_\Q(\Lambda_0)$ such that:
\begin{eqnarray*}
&&{\rm (G1)}\qquad G(\mu) \equiv v_\mu(\mod qL),\\
&&{\rm (G2)}\qquad \overline{G(\mu)} = G(\mu).
\end{eqnarray*}
\end{theorem}
The basis $\{G(\mu)\}$ is called the lower
global crystal basis of $V_\Q(\Lambda_0)$
(or simply of $V(\Lambda_0)$).

An efficient combinatorial means of computing the lower global
crystal basis of $V(\Lambda_0)$ was given by \cite{LLT}.
In general, for each $\mu\in\Pi_n(m)$,
this algorithm yields an expansion:
\begin{equation}\label{qExpand}
G(\mu)=\sum_{\lambda\in\Pi(m)} d_{\lambda\mu}(q)v_\lambda,
\end{equation}
where each $d_{\lambda\mu}(q)\in\Z[q]$.
In the case $n=2$, for partitions of weight 5, the $d_{\lambda\mu}$
are given in the following table, where the columns are indexed
by $\mu$ and the rows by $\lambda$.
\begin{equation}\label{LGCBmatrix}
\vcenter{\def\c{\cdot}\halign{\strut\tabskip=5mm $#$\hfil&&\hfil$#$\hfil\cr
&(5) &(41) &(32)\cr
\noalign{\vskip1mm}
(5)  &1&\c&\c\cr
(41) &\c&1&\c\cr
(32) &\c&\c&1\cr
(31^2) &q&\c&q\cr
(2^21) &\c&\c&q^2\cr
(21^3) &\c&q&\c\cr
(1^5)  &q^2&\c&\c\cr}}
\end{equation}
The similarity with (\ref{Hmatrix}) led the authors of \cite{LLT}
to conjecture the following result,
which was later proved by Ariki \cite{Ar2} and  by Grojnowski.

\begin{theorem}\label{LLTconjecture}
Let $d_{\lambda\mu}(q)$ be given by {\rm (\ref{qExpand})}.
Then the multiplicity of the irreducible module $D^\mu$
as a composition factor of the Specht module
$S^\lambda$ of $H_m(v)$ when $v$ is a primitive $n$th root of
unity is given by $d_{\lambda\mu}(1)$.
\end{theorem}

It is interesting to reformulate this theorem in terms of 
Grothendieck groups \cite{LLT}. Let $K_0(H_m(\sqrt[n]{1}))$ 
denote the Grothendieck group of the category of finitely 
generated projective $H_m(\sqrt[n]{1})$-modules. This is a 
free abelian group generated by the isomorphism classes of 
the indecomposable direct summands of $H_m(\sqrt[n]{1})$. 
If the $\slchap_n$-module $V(\Lambda_0)$ is regarded, via 
Robinson's $i$-induction and $i$-restriction operators, as 
the sum of Grothendieck groups 
$$
V(\Lambda_0) = {\cal K} = \bigoplus_m K_0(H_m(\sqrt[n]{1})),
$$
then the theorem states that the global lower crystal basis 
(at $q=1$) coincides with the canonical basis of ${\cal K}$ 
given by the classes of indecomposable summands of 
$H_m(\sqrt[n]{1})$.

Although only the values of $d_{\lambda\mu}(q)$ at $q=1$ are
required to determine the composition factors of $S^\lambda$,
it is conjectured in \cite{LLT}, that the coefficients of the 
polynomials $d_{\lambda\mu}(q)$ are positive and provide 
information on the order of the composition factors in the 
composition series of $S^\lambda$.  

\subsection{The upper global crystal basis}

The upper global crystal basis of $V(\Lambda_0)$ is defined as 
the basis adjoint to the lower global crystal basis with respect
to the so-called Shapovalov form of $V(\Lambda_0)$. This bilinear 
form is characterized by the following properties
$$
(v_\emptyset,v_\emptyset)=1,\qquad (q^hu,v)=(u,q^hv),
\qquad (f_iu,v)=(u,e_iv),
$$
for all $u,v\in V(\Lambda_0)$, $h\in P^\vee$, and $0\le i<n$.
Thus the upper global crystal basis is the unique basis 
$\{G^{up}(\mu)\mid \mu\in \Pi_n\}$ of $V(\Lambda_0)$ such that
$$
(G^{up}(\mu) , G(\nu) ) = \delta_{\mu \nu}, \quad (\mu,\nu \in \Pi_n).
$$
The interpretation of the upper basis in terms of Grothendieck 
groups of Hecke algebras is the following.
Let $G_0(H_m(\sqrt[n]{1}))$ be the Grothendieck group
of the category of finitely generated $H_m(\sqrt[n]{1})$-modules.
The elements of $G_0(H_m(\sqrt[n]{1}))$ are 
classes $[M]$ of modules, where 
$[M_1]=[M_2]$ if and only if
the composition factors of $M_1$
occur in $M_2$ with identical multiplicity.
(The order of the composition factors in the series is disregarded).
The sum is defined by $[M]+[N]=[M\oplus N]$. 
It is known that this is a free abelian group
with basis the set $\{[D^\mu]\}$ of classes of irreducible
$H_m(\sqrt[n]{1})$-modules.

Then the direct sum
$$
{\cal G}=\bigoplus_m G_0(H_m(\sqrt[n]{1}))
$$
endowed with Robinson's $i$-induction and $i$-restriction
operators becomes a $\slchap_n$-module isomorphic to
$V(\Lambda_0)$ \cite{LLT}. 
It has a canonical basis,
namely $\{[D^\mu],\ \mu \in \Pi_n\}$.
Another reformulation of Theorem~\ref{LLTconjecture}
is that this basis coincides with the upper global crystal
basis (at $q=1$).

In particular, we have the following 
interpretation of the complete restriction operator in
terms of the upper global basis of $V(\Lambda_0)$.

\begin{theorem}\label{GrothTheorem} 
Let $\lambda\in\Pi_n(m)$ and write
$$
\left( \sum_{i=0}^{n-1} e_i \right) G^{up}(\lambda)
= \sum_{\mu\in\Pi_n(m-1)} c_{\lambda\mu}(q) G^{up}(\mu).
$$
Then
$$
[D^\lambda\!\downarrow^{H_m}_{H_{m-1}}]
= \bigoplus_{\mu\in\Pi_n(m-1)} c_{\lambda\mu}(1) [D^\mu].
$$
\end{theorem}

\section{The Jantzen-Seitz problem for
	   Hecke algebras at roots of unity}

The results of \cite{LLT} and \cite{Ar2} provide us with a correspondence between 
the theory of modular representations of Hecke algebras and solvable 
lattice models. Under this correspondence, the Hecke JS problem is reformulated 
as a problem that is already solved in the context of exactly solvable 
lattice models, namely, the restriction of a solid-on-solid model to 
the simplest possible non-trivial restricted one. 

Let $\lambda \in \Pi_n$. One can write
$\lambda=(\lambda_1^{a_1},\lambda_2^{a_2},\ldots,\lambda_r^{a_r})$,
where $\lambda_1>\lambda_2>\cdots>\lambda_r>0$ and
$0<a_i<n$ for $i=1,2,\ldots,r$.
We define the set $JS(n)$ as the set of all $\lambda \in \Pi_n$ 
such that
$$
a_i+\lambda_i-\lambda_{i+1}+a_{i+1}=0 \pmod n,
$$
for $i=1,2,\ldots,r-1$.

\begin{theorem}\label{JStheorem}
Let $\lambda\in\Pi_n(m)$.
Then $D^\lambda\!\downarrow^{H_m}_{H_{m-1}}$ is irreducible
if and only if $\lambda\in JS(n)$.
\end{theorem}

\Proof
It is known
(see \cite{Ka3} 5.3.8.--5.3.10.) that
\begin{equation}\label{magic}
e_i\,G^{up}(\lambda) = [\varepsilon_i(\lambda)]_q
G^{up}(\tilde e_i \lambda) 
+ \sum_{\mu\ne\lambda} E_{\lambda\mu}^i \,G^{up}(\mu)
\end{equation}
where $E_{\lambda\mu}^i$ is a Laurent polynomial in $q$
invariant under $q\mapsto q^{-1}$, and such that:
\begin{displaymath}
\lim_{q\rightarrow 0} 
{E_{\lambda\mu}^i \over [\varepsilon_i(\lambda)]_q}=0\,.
\end{displaymath}
If $\varepsilon_i(\lambda)=1$ whereupon
$[\varepsilon_i(\lambda)]_q=1$, this implies that each
$E^i_{\lambda\nu}=0$ since for all $k$, the coefficients of $q^k$
and $q^{-k}$ in $E^i_{\lambda\mu}$ are equal.
Thereupon
$e_i\,G_q^{up}(\lambda) = G_q^{up}(\nu)$ where $\nu = \tilde e_i \lambda$.
On the other hand, if $e_i\,G_q^{up}(\lambda) = G_q^{up}(\nu)$
for some $\nu$ then since $[x]_q=q^{x-1}+q^{x-3}+\cdots+q^{-x+1}$,
(\ref{magic}) necessarily implies that
$\varepsilon_i(\lambda)=1$ whereupon $\nu = \tilde e_i \lambda$.

By Theorem \ref{GrothTheorem},
$D^\lambda\!\downarrow$ is irreducible if and only if
$\sum e_i G^{up}(\lambda) = G^{up}(\nu)$ for some $\nu$.
Therefore, by the above,
$D^\lambda\!\downarrow$ is irreducible if and only if
$\varepsilon_j(\lambda)=1$ for some $j\in\{0,1,\ldots ,n-1\}$
and $\varepsilon_i(\lambda) = 0$ for $i\ne j$.

Now consider the branching function $b_{\Lambda_j\,\Lambda_0}^\Lambda(q)$.
By Corollary \ref{TPcorollary}, the coefficient of
$q^k$ in this function is equal to the number of vertices $b$ 
of the crystal graph of $B(\Lambda_0)$ such that:
\begin{eqnarray*}
&&{\rm wt}(b) = \Lambda - \Lambda_j - k\delta\,; \\
&&\varepsilon_j(b)\le1\,;\\
&&\varepsilon_i(b)=0 \quad (i\not = j).
\end{eqnarray*}

It follows from this discussion that the partitions $\lambda$
such that $D^\lambda\downarrow$ is irreducible can be identified
with the set of vertices $b$ of the crystal graph $B(\Lambda_0)$
which contribute to a branching function
$b_{\Lambda_j\,\Lambda_0}^\Lambda(q)$ 
for some $j$ and some $\Lambda\in P_2^+$. By \ref{FOWTheorem} and 
\ref{ThJMMO}, this set of partitions is precisely $JS(n)$, which 
completes the proof.

\centerline{\epsfbox{sl3_0.eps}}

\begin{example}{\rm
The argument considered in the above proof may be illustrated by
reference to the $\emptyset$-connected component of the
crystal graph $\Gamma_3$ of $\slchap_3$ given
(up to weight 8) in Fig. 5.

\noindent
Here, those partitions $\lambda$ for which $\lambda\in JS(n)$
have been highlighted with an asterisk. As in the above proof, 
these partitions correspond to the nodes $b$ in this crystal 
graph for which, for some $j$, $\varepsilon_j(b)\le1$ and 
$\varepsilon_i(b)=0$ for all $i\ne j$.}
\end{example}

We will now give the interpretation in terms of
Hecke algebras (or symmetric groups) of the splitting
$$
JS(n) = \coprod_{j, k} FOW(n, j, k).
$$
To do this we recall the notions of $n$-core and
$n$-weight of a partition \cite{JK}.

If $\lambda$ is a partition, then a sequence of $n$ adjoining
nodes within $\lambda$ that begins at the rightmost end of a row,
passes to the node directly below if one exists, else passes to the node
directly on the left, and finishes at the bottom of a column,
is known as an $n$-rim-hook of $\lambda$.
For example, in the case $\lambda=(754^2)$, there are two
$5$-rim-hooks. In the following diagram, one comprises the nodes
labelled 3,4,5,6 and 7, and the other comprises the nodes labelled
5,6,7,8 and 9:
$$
\smyoungd{
\multispan{15}\hrulefill\cr
&&&&&&&&&2&&1&&0&\cr
\multispan{15}\hrulefill\cr
&&&&&&&4&&3&\cr
\multispan{11}\hrulefill\cr
&&&&&&&5&\cr
\multispan{9}\hrulefill\cr
&9&&8&&7&&6&\cr
\multispan{9}\hrulefill\cr}.
$$
(We also see that $\lambda$ has two $3$-rim-hooks and
four $4$-rim-hooks.)

 If $\lambda$ has no $n$-rim-hooks then $\lambda$ is termed
an $n$-{\it core}.

 Removing an $n$-rim-hook from $\lambda$ produces a valid partition.
Thus a process of successive $n$-rim-hook removals may be
carried out, the process eventually terminating with an $n$-core.
It may be shown that the $n$-core $\mu$ so obtained is independent
of the order of $n$-rim-hook removals (see \cite{JK}).
This unique $\mu$ is termed the {\it $n$-core of $\lambda$}.
For example, we find that $\lambda=(754^2)$ has $3$-core
$(421^2)$, $4$-core $\emptyset$ and $5$-core $(2^21)$.

The {\it $n$-weight} of $\lambda$ is the number of $n$-rim-hook
removals that are required to produce its $n$-core.
Thus $\lambda=(754^2)$ has 3-weight 4, 4-weight 5 and 5-weight~3.

The importance of the $n$-core derives from the
Nakayama \lq\lq Conjecture\rq\rq\ for $\SG_m$ (see \cite{JK})
and was extended to $H_m$ in \cite{DJ2}.

Let $JS(n,\mu,d)$  be the subset of $JS(n)$
comprising those partitions with $n$-core $\mu$ and $n$-weight $d$.
Define the generating series
$$
\chi_{n,\mu}(q)=\sum_{d\ge0} \#JS(n,\mu,d) q^d.
$$

\goodbreak
\begin{theorem}\label{JSCor}
\hfil
\begin{enumerate}
\item If $\lambda\in JS(n)$ then the $n$-core $\mu$ of $\lambda$
is a rectangular partition $\mu=(k^l)$ such that $k+l\le n$
(it is assumed here that if either $k=0$ or $l=0$ then
$(k^l)$ means the empty partition).
\item If $k\ne0\ne l$ then:
$$
\chi_{n,(k^l)}(q)= 
q^{-s}\,b_{\Lambda_{k-l}\,\Lambda_0}^{\Lambda_k + \Lambda_{-l}}(q) 
\,,
$$
where $s=\min(k,l)$.
\item
$$
\chi_{n,\emptyset}(q)=
\left(\,\sum_{k=0}^{n-1}
b_{\Lambda_k\,\Lambda_0}^{\Lambda_k + \Lambda_0}(q)\right)-(n-1)\,.
$$
\end{enumerate}
\end{theorem}

\Proof In the tensor product $V(\Lambda_j)\otimes V(\Lambda_0)$,
all highest weights are of the form $\Lambda_k+\Lambda_{-l}-e\delta$ 
for $k-l=j\mod n$ (for later convenience we use $-l$ and not $l$ here).
We take $0\le k,l<n$ here and can also assume that $k\le -l\mod n$
whereupon $l+k\le n$, and $l=0$ only if $k=0$.
By definition, the multiplicity of 
$V(\Lambda_k+\Lambda_{-l}-e\delta)$ is given by the coefficient of 
$q^e$ in $b_{\Lambda_j,\Lambda_0}^{\Lambda_k+\Lambda_{-l}}(q)$ and 
hence, by Theorem \ref{FOWTheorem} and Theorem \ref{ThJMMO}, by the 
number of partitions $\lambda\in JS(n)$ for which 
$\wt(\lambda)=\Lambda_k+\Lambda_{-l}-\Lambda_j-e\delta$.
We claim that the $n$-core of such a $\lambda$ is the rectangular
partition $\mu=(k^l)$, for which, by Lemma \ref{EwtLemma}, we calculate
$\wt(\mu)=\Lambda_k+\Lambda_{-l}-\Lambda_{k-l}-s\delta=
\Lambda_k+\Lambda_{-l}-\Lambda_j-s\delta$,
where $s=\min(k,l)$ is the multiplicity of the colour charge 0 in $\mu$.
To establish this claim first note that for every string of weights
$\Lambda,\Lambda-\delta,\Lambda-2\delta,\ldots,$ of the
$\slchap_n$-module $V(\Lambda_0)$ (where $\Lambda+\delta$ is not
a weight of $V(\Lambda_0)$), those partitions having these
weights have the same $n$-core. This $n$-core has weight $\Lambda$.
Then since $\mu=(k^l)$ with $k+l\le n$, is manifestly an $n$-core,
it follows that it is the $n$-core of $\lambda$, thence proving
part 1.

Part 2 follows immediately since in the case $k\ne0$ (so that
$l\ne0$), the partitions that enumerate the branching function
$b_{\Lambda_{k-l}\Lambda_0}^{\Lambda_k+\Lambda_{-l}}$
are precisely those elements $\lambda\in JS(n)$ having weight 
$\wt(\lambda)=\Lambda_k+\Lambda_{-l}-\Lambda_{k-l}-e\delta$,
for some $e$, and hence $n$-core $(k^l)$.

For $k=0$ and arbitrary $l$, each partition enumerating
$b_{\Lambda_{-l}\Lambda_{0}}^{\Lambda_{0}+\Lambda_{-l}}$
has $n$-core $\emptyset$, and hence contributes
to $\chi_{n,\emptyset}$.
However, the empty partition occurs for each $l$, hence an adjustment
of $n-1$ is needed after summing over all $l$.
No other partition $\lambda$ is repeated since, as indicated by
Theorem \ref{FOWTheorem},
the $b_{\Lambda_{-l}\Lambda_{0}}^{\Lambda_{0}+\Lambda_{-l}}$
to which it contributes is uniquely determined by
$-l\mod n=(\lambda_1-a_1)\mod n$.
(The summation over $-l$ is replaced by one over $k$ to give
the final result).

\begin{example}{\rm

To illustrate this result, consider again the case $n=3$,
where we have the following branching functions (to three terms):
\begin{equation}
\begin{array}{lll}
b_{\Lambda_0,\Lambda_0}^{2\Lambda_0}&=&1+q^2+\cdots;\\[1mm]
b_{\Lambda_0,\Lambda_0}^{\Lambda_1+\Lambda_2}&=&q+2q^2+2q^3+\cdots;\\[1mm]
b_{\Lambda_1,\Lambda_0}^{\Lambda_1+\Lambda_0}&=&1+q+2q^2+\cdots;\\[1mm]
b_{\Lambda_1,\Lambda_0}^{2\Lambda_2}&=&q+q^2+2q^3+\cdots;\\[1mm]
b_{\Lambda_2,\Lambda_0}^{\Lambda_2+\Lambda_0}&=&1+q+2q^2+\cdots;\\[1mm]
b_{\Lambda_2,\Lambda_0}^{2\Lambda_1}&=&q+q^2+2q^3+\cdots.
\end{array}
\end{equation}
These are calculated from Theorem \ref{FOWTheorem} which leads to 
the enumeration of the nodes of Fig.~5. labelled by asterisks. The 
only rectangular $3$-cores are $\emptyset$, $(1)$, $(2)$ and $(1^2)$.
Using Theorem \ref{JSCor}, we thus obtain:
\begin{equation}
\begin{array}{lll}
\chi_{n,\emptyset}&=&1+2q+5q^2+\cdots;\cr
\chi_{n,(1)}&=&1+2q+2q^2+\cdots;\cr
\chi_{n,(2)}&=&1+q+2q^2+\cdots;\cr
\chi_{n,(1^2)}&=&1+q+2q^2+\cdots.\cr
\end{array}
\end{equation}
\noindent These correspond to the following sets:
\begin{equation}
\begin{array}{lll}
JS(3,\emptyset,0)&=&\{\emptyset\};\cr
JS(3,\emptyset,1)&=&\{(3),(21)\};\cr
JS(3,\emptyset,2)&=&\{(6),(51),(3^2),(41^2),(321)\},\cr
JS(3,(1),0)&=&\{(1)\};\cr
JS(3,(1),1)&=&\{(4),(2^2)\};\cr
JS(3,(1),2)&=&\{(7),(43)\},\cr
JS(3,(2),0)&=&\{(2)\};\cr
JS(3,(2),1)&=&\{(5)\};\cr
JS(3,(2),2)&=&\{(8),(3^21^2)\},\cr
JS(3,(1^2),0)&=&\{(1^2)\};\cr
JS(3,(1^2),1)&=&\{(32)\};\cr
JS(3,(1^2),2)&=&\{(62),(44)\}.\cr
\end{array}
\end{equation}}
\end{example}

\newpage

\section{Discussion}

We have explained why the same set of combinatorial objects
arises in two different contexts, namely the modular
representations of symmetric groups and Hecke algebras,
and the RSOS solvable lattice models based on the
coset algebras $C[n, 1, 1]$. 
In fact the Jantzen-Seitz problem and the problem of 
evaluating one-dimensional configurations associated with
these models are both equivalent to the computation
of branching functions of tensor products of level 1
representations of $\slchap_n$.

It is very natural to ask what happens when one considers
tensor products of higher level representations.
By Ariki's theorem \cite{Ar2}, the level $\ell$ representations of
$\slchap_n$ can be interpreted as Grothendieck groups of
some cyclotomic Hecke algebras \cite{AK,BM}.
In particular for level 2 one gets Hecke algebras
of type $B$ at roots of unity.  
On the other hand, removing the level 1 restriction, one
obtains the more general RSOS models associated with
cosets $C[n, \ell, 1]$.
This generalized correspondence is currently under investigation
and we hope to report on this subject soon.

\footnotesize
\section*{Acknowledgements}

We wish to thank Professor Christine Bessenrodt for discussions that 
initiated this work, and Dr Ole Warnaar for an earlier collaboration 
on which it is partly based. This work was financially supported by 
the Australian Research Council.

{\noindent \it Note added --- }
A forthcoming preprint, \cite{BO}, contains (among other things) an 
elementary purely combinatorial proof of part $(i)$ of \ref{JSCor}.

\normalsize

\begin{appendix}

\section{The Specht module $S^\lambda$.}


\def\hbf#1{{\hbox{\scriptsize\bf #1}}}

\def\exyoungd#1#2#3#4#5#6{%
 \smyoungd{\multispan{13}\hrulefill\cr
 &1&&#1&&#2&&10&&4&&12&\cr
 \multispan{13}\hrulefill\cr
 &6&&#3&&#4&\cr
 \multispan{7}\hrulefill\cr
 &9&&#5&&#6&\cr
 \multispan{7}\hrulefill\cr
 &13&\cr
 \multispan{3}\hrulefill\cr}
}

\def\eyyoungd#1#2#3#4#5{%
 \smyoungd{\multispan{7}\hrulefill\cr
 &#1&&#2&&#3&\cr
 \multispan{7}\hrulefill\cr
 &#4&&#5&\cr
 \multispan{5}\hrulefill\cr}
}

In this Appendix we give an explanation of the construction of
the Specht module representations of $H_m(v)$ based on the
account given in \cite{KWy} which makes use of Young tableaux.
An alternative construction is given in \cite{DKLLST}.

Let $\lambda\in\Pi(m)$. 
Filling the $m$ nodes of $\lambda$ with elements of $\{1,2,\ldots,m\}$
so that no element appears more than once,
yields what is known as a Young tableau.
A Young tableau of shape $\lambda$ is typically denoted $t^\lambda$.

A Young tableau $t^\lambda$ for which the entries increase
from left to right across each row is said to be row-standard.
A Young tableau $t^\lambda$ for which the entries increase
from top to bottom down each column is said to be column-standard.

A Young tableau $t^\lambda$ which is both row-standard
and column-standard is said to be standard.

The standard Young tableau of shape $\lambda$ for which the
entries $1,2,\ldots,m$
occur successively down first the leftmost column then down subsequent
columns taken in turn from left to right is denoted $t^\lambda_-$.

Let $1\le a,b\le m$. The entry $a$ is said to precede $b$
in the tableau $t^\lambda$ if $a$ precedes $b$ on reading
the entries of $t^\lambda$ down first the leftmost column of
$t^\lambda$ and then down the remaining columns taken in turn from
left to right.
($t^\lambda_-$ is such that $a$ precedes $b$ if and only if $a<b$.)

\begin{example} {\rm
In the case when $\lambda=(42)$, there are of
course $6!=720$ Young tableaux. Of these, 15 are row-standard,
180 are column standard and just the following nine are standard:
$$
\smyoungd{ \multispan9\hrulefill\cr
&1&&3&&5&&6&\cr \multispan9\hrulefill\cr
&2&&4&\cr \multispan5\hrulefill\cr}, \quad
\smyoungd{ \multispan9\hrulefill\cr
&1&&3&&4&&6&\cr \multispan9\hrulefill\cr
&2&&5&\cr \multispan5\hrulefill\cr}, \quad
\smyoungd{ \multispan9\hrulefill\cr
&1&&2&&5&&6&\cr \multispan9\hrulefill\cr
&3&&4&\cr \multispan5\hrulefill\cr}, \quad
\smyoungd{ \multispan9\hrulefill\cr
&1&&2&&4&&6&\cr \multispan9\hrulefill\cr
&3&&5&\cr \multispan5\hrulefill\cr}, \quad
\smyoungd{ \multispan9\hrulefill\cr
&1&&3&&4&&5&\cr \multispan9\hrulefill\cr
&2&&6&\cr \multispan5\hrulefill\cr},
$$ $$
\smyoungd{ \multispan9\hrulefill\cr
&1&&2&&3&&6&\cr \multispan9\hrulefill\cr
&4&&5&\cr \multispan5\hrulefill\cr}, \quad
\smyoungd{ \multispan9\hrulefill\cr
&1&&2&&4&&5&\cr \multispan9\hrulefill\cr
&3&&6&\cr \multispan5\hrulefill\cr}, \quad
\smyoungd{ \multispan9\hrulefill\cr
&1&&2&&3&&5&\cr \multispan9\hrulefill\cr
&4&&6&\cr \multispan5\hrulefill\cr}, \quad
\smyoungd{ \multispan9\hrulefill\cr
&1&&2&&3&&4&\cr \multispan9\hrulefill\cr
&5&&6&\cr \multispan5\hrulefill\cr}.
$$
The first tableau listed here is $t^{(42)}_-$.
}
\end{example}

\begin{definition}[The Specht module]{\rm \cite{KWy}}
The Specht module $S^\lambda$ of $H_m(v)$
is spanned by vectors $u_{t^\lambda}$ indexed by Young tableaux
of shape $\lambda$ subject to the column and Garnir relations
(defined below), and on which the action of $H_m(v)$ is generated by:
$$
T_i u_{t^\lambda} =
\cases{
u_{x^\lambda}
&if $i$ precedes $i+1$ in $t^\lambda$;\cr
v u_{x^\lambda} + (v-1) u_{t^\lambda}
&if $i+1$ precedes $i$ in $t^\lambda$,\cr}
$$
where the tableau $x^\lambda$ is obtained from $t^\lambda$ by
interchanging $i$ and $i+1$.
\end{definition}

\noindent Note that the Specht modules for the symmetric group $\SG_m$
are obtained simply by setting $v=1$ throughout this Appendix.
In particular, we see that the action of $\SG_m$
on $u_{t^\lambda}$ does not depend on whether $i$ precedes $i+1$.
Then the action of the generator $s_i$ of $\SG_m$ simply interchanges
$i$ and $i+1$ in all cases.

\begin{definition}[Column relations]
If $x^\lambda$ differs from $z^\lambda$ only in that a single pair of
entries within a column are transposed then:
$$
u_{z^\lambda}=-v_{x^\lambda}.
$$
\end{definition}

\setsmyoungsize{11pt}{9pt}

\begin{example}{\rm
For example (denoting $u_{t^\lambda}$ by $t^\lambda$ for
typographical reasons),
$$
\exyoungd{8}{5}{11}{3}{2}{7}
=-\:\exyoungd{8}{5}{2}{3}{11}{7}
=-\:\exyoungd{2}{3}{8}{5}{11}{7}.
$$
}
\end{example}

Let $z^\lambda$ be column-standard but not row-standard.
Then an adjacent pair of entries exists with that on
the left greater than that on the right. Consider these two entries
together with all those below the left one and all those above
the right one (see the example below).
Now form all possible tableaux $t^\lambda$ by permuting these entries
in all ways such that the permuted entries are increasing
down the portions of each of the two columns being considered.

\begin{definition}[Garnir relations]
The Garnir relation is  the following expression in which the sum
is over all such tableaux:
$$
(-v)^{l(w_{z^\lambda})}
\sum_{t^\lambda} (-v)^{-l(w_{t^\lambda})}
u_{t^\lambda}=0,
$$
where $w_{t^\lambda}\in \SG_m$ maps $t^\lambda_-$ to $t^\lambda$
and $l(w)$ is the length of the permutation $w\in\SG_m$.
(The factor at the front is included to make all coefficients
postive powers of $v$ and the coefficient of $u_{z^\lambda}$ one).
\end{definition}

\noindent This enables $u_{z^\lambda}$ to be written in terms of other
tableaux.

\begin{example}{\rm
Consider the following tableau:
$$
z^\lambda=
\exyoungd{2}{\hbf3}{\hbf8}{\hbf5}{\hbf{11}}{7}.
$$
If we consider the pair of entries 8 and 5 which violate the
row-standard condition, we are led to permuting the entries which
are highlighted here.
The following Garnir relation results:
\begin{displaymath}
\begin{array}{l}
\exyoungd{2}{\hbf3}{\hbf8}{\hbf5}{\hbf{11}}{7}
-v\:\exyoungd{2}{\hbf3}{\hbf5}{\hbf8}{\hbf{11}}{7}
+v^2\:\exyoungd{2}{\hbf5}{\hbf3}{\hbf8}{\hbf{11}}{7}\\[7mm]
\qquad\qquad
+v^2\:\exyoungd{2}{\hbf3}{\hbf5}{\hbf{11}}{\hbf8}{7}
-v^3\:\exyoungd{2}{\hbf5}{\hbf3}{\hbf{11}}{\hbf8}{7}
+v^4\:\exyoungd{2}{\hbf8}{\hbf3}{\hbf{11}}{\hbf5}{7}
=0.\\
\end{array}
\end{displaymath}
}
\end{example}

\noindent
By employing a lexicographical ordering on the set of tableaux,
it may be seen that each element of $S^\lambda$ may be written
as a linear combination of elements $u_{t^\lambda}$, by repeated
use of the column and Garnir relations. In fact, we have the
following:

\begin{theorem}{\rm \cite{DJ1}}
\hfil
\begin{enumerate}
\item The set
$$
\{u_{t^\lambda}\mid t^\lambda {\rm\ is\ a\ standard\ tableau}\}
$$
is a basis for $S^\lambda$.
\item If $v$ is not a root of unity, then
$$
\{S^\lambda\mid\lambda\in\Pi(m)\}
$$
is a complete set of mutually inequivalent irreducible
$H_m(v)$-modules.
\end{enumerate}
\end{theorem}

\noindent Note that the construction implies that each matrix
element of the resulting representation matrices is an
element of $\Z[v]$. This implies that these representation
matrices remain well defined on specialising $v$ to any value
or on taking the entries modulo some field characteristic,
and indeed provide representation matrices in these circumstances.

\setsmyoungsize{9pt}{9pt}

\begin{example}{\rm
To illustrate the construction of explicit representation matrices
using the Specht module approach,
consider representing $T_1\in H_5(v)$ in the Specht
module $S^{(32)}$, by acting on each $u_{t^{(32)}}$ for which $t^{(32)}$
is standard (once more $u_{t^\lambda}$ is written as $t^\lambda$):
$$
\begin{array}{l}
{ T_1} \:\eyyoungd13524
= \eyyoungd23514
= -\:\eyyoungd13524,\\[3mm]
{ T_1} \:\eyyoungd12534
= \eyyoungd21534
= v\:\eyyoungd12534-v^2\:\eyyoungd13524,\\[3mm]
{ T_1} \:\eyyoungd13425
= \eyyoungd23415
= -\:\eyyoungd13425,\\[3mm]
{ T_1} \:\eyyoungd12435
= \eyyoungd21435
= v\:\eyyoungd12435-v^2\:\eyyoungd13425,\\[3mm]
{ T_1} \:\eyyoungd12345
= \eyyoungd21345
= v\:\eyyoungd12345-v^2\:\eyyoungd14325
=v \:\eyyoungd12345
-v^3 \:\eyyoungd13425
+v^4 \:\eyyoungd13524.\\
\end{array}
$$
Here, column relations have been used in the first and third calculations,
and Garnir relations have been used in the second, fourth and last (twice),
to express each result in terms of the standard tableaux.
Consequently, in the representation labelled by the partition $(32)$,
$T_1$ is mapped to the matrix (where zeros are denoted by dots):
\begin{displaymath}
\pmatrix{
-1&-v^2&.&.&v^4\cr
.&v&.&.&.\cr
.&.&-1&-v^2&-v^3\cr
.&.&.&v&.\cr
.&.&.&.&v\cr}.
\end{displaymath}
The matrices representing the generators $T_i$ of $H_m(v)$ in each
irreducible representation for $m\le5$ given
in \cite{KWy} have been produced in a similar way.
}
\end{example}

\end{appendix}


\newpage \small
\begin{thebibliography}{99}


\bibitem{ABF}{\sc G.E.~Andrews, R.J.~Baxter} and {\sc P.J.~Forrester},
{\it Eight-vertex SOS model and generalized Rogers-Ramanujan-type 
     identities},
J. Stat. Phys. {\bf 35} (1984), 193--266.



\bibitem{Ar2}{\sc S.~Ariki},
{\it On the decomposition numbers of the Hecke algebra of $G(m,1,n)$},
preprint, 1996.


\bibitem{AK}{\sc S.~Ariki} and {\sc K.~Koike},
{\it A Hecke algebra of $(\Z/rZ)\wr \SG_n$ and construction
of its irreducible representations},
Adv. Math. {\bf 106} (1994), 216-243.


\bibitem{B}{\sc R.J.~Baxter},
{\it Hard hexagons: exact solution},
J. Phys. {\bf A~13} (1980), L61--L70.


\bibitem{BaxterBook}{\sc R.J.~Baxter},
{\it Exactly solved models in statistical mechanics},
(1984) Academic Press. 

\bibitem{BPZ} {\sc A. A.~Belavin, A. M.~Polyakov} and 
              {\sc A. B.~Zamolodchikov},
{\it Infinite conformal symmetry in two-dimensional
quantum field theory},
Nucl. Phys. {\bf B~241} (1984) 333-380.


\bibitem{Be}{\sc D.~Benson},
{\it Some remarks on the decomposition numbers of the symmetric 
groups},
Proc. Symp. Pure Math. {\bf 47} (1987), 381-394.

\bibitem{BM}{\sc M.~Brou\'e} and {\sc G.~Malle},
{\it Zyklotomische Heckealgebren},
Asterisque {\bf 212} (1993), 119--189.


\bibitem{ChariPressleyBook} {\sc V.~Chari} and {\sc A.~Pressley},
{\it A guide to quantum groups},
Cambridge University Press, 1994.

\bibitem{DJKMOl} {\sc E.~Date, M.~Jimbo, A.~Kuniba, T.~Miwa} and 
                 {\sc M.~Okado},
{\it A new realization of the basic representation of $A_n^{(1)}$},
Lett. Math. Phys. {\bf 17} (1989), 51--54.

\bibitem{DJKMO} {\sc E.~Date, M.~Jimbo, A.~Kuniba, T.~Miwa} and 
                {\sc M.~Okado},
{\it Paths, Maya diagrams, and representations of $\slchap(r,C)$},
Adv. Stud. Pure Math. {\bf 19} (1989), 149--191.

\bibitem{DJMO1} {\sc E.~Date, M.~Jimbo, T.~Miwa} and 
               {\sc M.~Okado},
{\it Automorphic properties of local height probabilities for
integrable solid-on-solid models},
Phys. Rev. {\bf B 35} (1987), 2105--2107.

\bibitem{DJMO2}{\sc E.~Date, M.~Jimbo, T.~Miwa} and 
               {\sc M.~Okado},
{\it Solvable lattice models},
Proc. Symp. Pure Math. {\bf 49} (1989), 295--331.


\bibitem{DFJMN}{\sc B.~Davies, O.~Foda, M.~Jimbo, T.~Miwa} and
               {\sc A. Nakayashiki},
{\it Diagonalization of the $XXZ$ Hamiltonian by vertex operators},
Commun. Math. Phys. {\bf 151} (1993), 89--153.

\bibitem{DJ1} {\sc R.~Dipper} and {\sc G.D.~James},
{\it Representations of Hecke algebras of general linear groups},
Proc. London Math. Soc. {\bf 52} (1986), 20--52.

\bibitem{DJ2} {\sc R.~Dipper} and {\sc G.D.~James},
{\it Blocks and idempotents of Hecke algebras of general linear groups},
Proc. London Math. Soc. {\bf 54} (1987), 57--82.

\bibitem{DKLLST}{\sc G.~Duchamp, D.~Krob, A.~Lascoux, B.~Leclerc, T.~Scharf} 
            and {\sc J.-Y.~Thibon},
{\it Euler-Poincar\'e characteristic and polynomial
representations of Iwahori-Hecke algebras},
Publ. RIMS Kyoto. Univ. {\bf 31} (1995), 179--201.

\bibitem{FLOTW} {\sc O.~Foda, B.~Leclerc, M.~Okado, J.-Y.~Thibon}
            and {\sc T.A.~Welsh}
{\it Restricted solid-on-solid models, and Jantzen-Seitz
representations of Hecke algebras at roots of unity}, preprint, 
1997. 

\bibitem{FOW}{\sc O.~Foda, M.~Okado} and {\sc S.O.~Warnaar},
{\it A proof of polynomial identities of type 
$\widehat{sl}(n)_1\otimes\widehat{sl}(n)_1/\widehat{sl}(n)_2$},
J. Math. Phys. {\bf 37} (1996), 965--986.

\bibitem{GKO}{\sc P.~Goddard, A.~Kent, D.~Olive}, 
{\it Virasoro algebras and coset space models},
Phys. Lett. {\bf B~152} (1985), 88-93.

\bibitem{Gr}{\sc I.~Grojnowski}, 
{\it Representations of affine Hecke algebras (and affine
quantum $GL_n$) at roots of unity},
Internat. Math. Research Notices, {\bf 5} (1994), 215-217.

\bibitem{Hay}{\sc T.~Hayashi},
{\it $q$-analogues of Clifford and Weyl algebras - spinor and
oscillator representations of quantum enveloping algebras},
Commun. Math. Phys. {\bf 127} (1990), 129--144.

\bibitem{ISZBook}{\sc C.~Itzykson, H.~Saleur} and {\sc J.-B.~Zuber},
{\it Conformal invariance and applications to statistical mechanics},
World Scientific, 1988.

\bibitem{Ja}{\sc G.D.~James},
{\it The representation theory of the symmetric groups},
Lecture Notes in Mathematics {\bf 682},
Springer, 1978.


\bibitem{JK}{\sc G.D.~James} and {\sc A.~Kerber},
{\it The representation theory of the symmetric group},
Addison-Wesley, 1981.


\bibitem{JS}{\sc J.C.~Jantzen} and {\sc G.M.~Seitz},
{\it On the representation theory of the symmetric groups},
Proc. London Math. Soc. {\bf 65} (1992), 475--504.

\bibitem{JMMN}{\sc M.~Jimbo, K.~Miki, T.~Miwa} and 
              {\sc A. Nakayashiki},
{\it Correlation functions of the {X}{X}{Z} model for ${\Delta}<-1$},
Phys. Lett. {\bf A~168} (1992), 256--263.

\bibitem{JMMO}{\sc M.~Jimbo, K.~Misra, T.~Miwa} and 
              {\sc M.~Okado},
{\it Combinatorics of representations of $U_q(\slchap(n))$ at $q=0$},
Commun. Math. Phys. {\bf 136} (1991), 543--566.

\bibitem{JM}{\sc M.~Jimbo} and {\sc T.~Miwa},
{\it Solitons and infinite dimensional Lie algebras},
Publ. RIMS Kyoto Univ. {\bf 19} (1983), 943--1001.

\bibitem{JimboMiwaBook}{\sc M.~Jimbo} and {\sc T.~Miwa},
{\it Algebraic analysis of solvable lattice models},
Providence, 1995.

\bibitem{JMO}{\sc M.~Jimbo, T.~Miwa} and {\sc M.~Okado},
{\it Solvable lattice models whose states are dominant integral
     weights of $A_{n-1}$}
Lett. Math. Phys. {\bf 14} (1987) 123--131.

\bibitem{KacBook}{\sc V.~Kac},
{\it Infinite dimensional Lie algebras}, 
3rd. ed., Cambridge, 1990.

\bibitem{KacRainaBook}{\sc V.~Kac} and {\sc A.K.~Raina},
{\it Bombay lectures on highest weight representations of
infinite dimensional Lie algebras},
World Scientific, 1987.

%

\bibitem{Ka1}{\sc M.~Kashiwara},
{\it Crystalizing the $q$-analogue of universal enveloping algebras},
Commun. Math. Phys. {\bf 133} (1990), 249--260.

\bibitem{Ka2}{\sc M.~Kashiwara},
{\it On crystal bases of the $q$-analogue of universal enveloping algebras},
Duke Math. J. {\bf 63} (1991), 465--516.

\bibitem{Ka3}{\sc M.~Kashiwara},
{\it Global crystal bases of quantum groups},
Duke Math. J. {\bf 69} (1993), 455--485.

\bibitem{KassBook}{\sc S.~Kass, R.V.~Moody, J.~Patera} and {\sc R.~Slansky},
{\it Affine Lie algebras, weight multiplicities, and branching rules},
University of California Press, 1990.

\bibitem{KWy}{\sc R.C.~King} and {\sc B.G.~Wybourne},
{\it Representations and traces of the Hecke algebras $H_n(q)$
of type $A_{n-1}$},
J. Math. Phys. {\bf33} (1992), 4--14.

\bibitem{Kl1}{\sc A.S. Kleshchev},
{\it On restrictions of irreducible modular representations
of semisimple algebraic groups and symmetric groups to natural 
subgroups I},
Proc. London Math. Soc. {\bf 69} (1994), 515-540.

\bibitem{Kl3}{\sc A.S. Kleshchev},
{\it Branching rules for the modular representations of symmetric groups III; 
some corollaries and a problem of Mullineux},
J. London Math. Soc. {\bf 2} (1995).

\bibitem{LLT1}{\sc A. Lascoux, B. Leclerc} and {\sc J.-Y. Thibon},
{\it Une conjecture pour le calcul des matrices de d\'ecomposition
des alg\`ebres de Hecke de type $A$ aux racines de l'unit\'e},
C. R. Acad. Sci. Paris {\bf 321} (1995), 511-516.

\bibitem{LLT}{\sc A.~Lascoux, B.~Leclerc} and {\sc J.-Y.~Thibon},
{\it Hecke algebras at roots of unity and crystal bases of
quantum affine algebras},
Commun. Math. Phys. (to appear).

\bibitem{LT}{\sc B.~Leclerc} and {\sc J.-Y.~Thibon},
{\it Canonical bases of $q$-deformed Fock spaces},
Internat. Math. Research Notices, {\bf 9} (1996), 447-456.

\bibitem{MartinBook}{\sc P.~Martin}
{\it Potts models and related problems in statistical mechanics},
(1991) World Scientific.

\bibitem{MM}{\sc K.C.~Misra} and {\sc T.~Miwa},
{\it Crystal base of the basic representation of $U_q(\widehat{\Sl}_n)$},
Commun. Math. Phys. {\bf 134} (1990), 79--88.

\bibitem{O} {\sc L.~Onsager},
{\it Crystal statistics I. A two dimensional model with an
order-disorder transition},
Phys. Rev. {\bf A~65} (1944), 117--149.

\bibitem{RC} {\sc A.~Rocha-Caridi},
{\it Vacuum vector representations
of the Virasoro algebra}, in ``Vertex Operators in Mathematics
and Physics'', 
MSRI Publication {\bf \# 3} 
(Springer, Heidelberg 1984), 451--473.

\bibitem{BO}{\sc Ch. Bessenrodt} and {\sc J. Olsson},
{\it Residue symbols and Jantzen-Seitz partitions},
to appear.

\end{thebibliography}
\end{document}